# SEMIPARAMETRIC ESTIMATION FOR STATIONARY PROCESSES WHOSE SPECTRA HAVE AN UNKNOWN POLE[1]

BY JAVIER HIDALGO

*London School of Economics*

We consider the estimation of the location of the pole and memory parameter, $\lambda^0$ and $\alpha$, respectively, of covariance stationary linear processes whose spectral density function $f(\lambda)$ satisfies $f(\lambda) \sim C|\lambda - \lambda^0|^{-\alpha}$ in a neighborhood of $\lambda^0$. We define a consistent estimator of $\lambda^0$ and derive its limit distribution $Z_{\lambda^0}$. As in related optimization problems, when the true parameter value can lie on the boundary of the parameter space, we show that $Z_{\lambda^0}$ is distributed as a normal random variable when $\lambda^0 \in (0, \pi)$, whereas for $\lambda^0 = 0$ or $\pi$, $Z_{\lambda^0}$ is a mixture of discrete and continuous random variables with weights equal to $1/2$. More specifically, when $\lambda^0 = 0$, $Z_{\lambda^0}$ is distributed as a normal random variable truncated at zero. Moreover, we describe and examine a two-step estimator of the memory parameter $\alpha$, showing that neither its limit distribution nor its rate of convergence is affected by the estimation of $\lambda^0$. Thus, we reinforce and extend previous results with respect to the estimation of $\alpha$ when $\lambda^0$ is assumed to be known a priori. A small Monte Carlo study is included to illustrate the finite sample performance of our estimators.

**1. Introduction.** Given a covariance stationary process $\{x_t\}$ observed at times $t = 1, 2, \ldots, n$, the search for cyclical components and their estimation and testing are of undoubted interest. This is motivated by the observed periodic behavior exhibited in many time series and manifested by sharp peaks of the spectral density estimate.

A well-known model capable of generating such a periodic behavior is the regression model

$$(1.1) \qquad x_t = \mu + \rho_1 \cos(\lambda^0 t) + \rho_2 \sin(\lambda^0 t) + \varepsilon_t,$$

Received March 1999; revised July 2004.

[1] Supported by the Economic and Social Research Council (ESRC) reference number R000238212.

where $\rho_1$ and $\rho_2$ are zero-mean uncorrelated random variables with the same variance and $\{\varepsilon_t\}$ is a stationary sequence of random variables independent of $\rho_1$ and $\rho_2$. Model (1.1) has enjoyed extensive use and different techniques have been proposed for the estimation of the frequency, amplitude and phase (see [6, 18, 19, 20, 35]). Extensions to a model with more than one periodic component have been examined by Quinn [29] and Kavalieris and Hannan [25], whose interest was also in testing the number of sinusoidal/cosinusoidal components. See also [30].

A second statistical model capable of exhibiting peaks in its spectral density function is the autoregressive AR(2) process

$$(1.2) \qquad (1 - a_1 L - a_2 L^2) x_t = \varepsilon_t$$

when the zeros of the polynomial $(1 - a_1 L - a_2 L^2)$ are complex, with $\lambda^0$ identified as $\arccos(\frac{a_1}{2\sqrt{-a_2}})$. Models (1.1) and (1.2) represent two extreme situations explaining cyclical behavior of the data and the peakedness of the spectral density function. Model (1.2) possesses a continuous spectral density function whereas model (1.1) has a spectral distribution function with a jump at the frequency $\lambda^0$. The cyclical component of the data remains constant or invariant with time in model (1.1), whereas the cyclical pattern of model (1.2) fades out with time fairly quickly.

Between these two extreme situations there exists a class of intermediate models in which the spectral density function of $x_t$ exhibits a pole at the frequency $\lambda^0$. For that purpose, define the spectral density function of $x_t$ as the function $f(\lambda)$ which satisfies the relationship

$$(1.3) \quad \gamma(j) = \mathrm{Cov}(x_t, x_{t+j}) = \int_{-\pi}^{\pi} f(\lambda) \cos(j\lambda)\, d\lambda, \qquad j = 0, 1, 2, \ldots.$$

We say that $f(\lambda)$ has a pole at $\lambda^0$ if

$$(1.4) \qquad f(\lambda) \sim C|\lambda - \lambda^0|^{-\alpha} \qquad \text{as } \lambda \to \lambda^0,$$

where $C \in (0, \infty)$, $\alpha \in (0, 1)$ is the memory parameter and "$\sim$" means that the ratio of the left- and right-hand sides tends to 1. One of the main objectives of this paper is the estimation of $\lambda^0$.

One model capable of generating such a cyclical behavior in the data has been proposed by Andel [2] and Gray, Zhang and Woodward [17] and defined as

$$(1.5) \qquad (1 - 2(\cos \lambda^0) L + L^2)^d x_t = \varepsilon_t,$$

where $L$ is the backshift operator, $d = \alpha/2$ for $\lambda^0 \in (0, \pi)$, whereas for $\lambda^0 = 0$ or $\pi$, $d = \alpha/4$. The model (1.5) was coined the Gegenbauer model by Gray, Zhang and Woodward [17], who extended it to the GARMA model where the innovations $\{\varepsilon_t\}$ follow an autoregressive moving average (ARMA) process, and it was later extended by Giraitis and Leipus [13] allowing for more than one pole or cyclical component. The GARMA process is characterized by having the spectral density function

$$(1.6) \quad f(\lambda) = \frac{\sigma^2}{2\pi}|1 - 2(\cos\lambda^0)e^{i\lambda} + e^{i2\lambda}|^{-2d}\left|\frac{a(e^{i\lambda};\theta)}{b(e^{i\lambda};\theta)}\right|^2, \qquad -\pi < \lambda \leq \pi,$$

where $\sigma^2 > 0$, and $a(\cdot)$ and $b(\cdot)$ are polynomials of finite degree, all of whose zeroes lie outside the unit circle. When $\lambda^0 = 0$, we have the more familiar FARIMA model, apparently originated by Adenstedt [1], and studied by Granger and Joyeux [16] and Hosking [23]. GARMA models are characterized by a stronger and more persistent cyclical behavior than ARMA models, that is, (1.2), but unlike model (1.1), their amplitude does not remain constant over time.

If the location of the pole $\lambda^0$ is known, then under some regularity conditions and a correct specification of the model, Whittle estimates of the parameters $\alpha$ (or $d$), $\theta$ and $\sigma^2$, for example model (1.6), are known to be $n^{1/2}$-consistent and asymptotically normal. In the case of Gaussianity or linear processes, this was shown by Fox and Taqqu [11], Dahlhaus [7] and Giraitis and Surgailis [15] when $\lambda^0 = 0$ and generalized by Giraitis and Leipus [13] and Hosoya [24] for $\lambda^0$ different from 0. All these papers assume that $f(\lambda)$ is fully specified by a set of parameters $(\alpha, \theta', \sigma^2)'$.

Although knowledge of $\lambda^0$ can be realistic in some time series data, with nonseasonal data that knowledge of $\lambda^0$ is not so clear. An example of the latter is when the practitioner is interested in estimating cycles in macroeconomic or geophysical data. Recently, Giraitis, Hidalgo and Robinson [12] have shown that Whittle estimates of $(\alpha, \theta', \sigma^2)'$ are asymptotically the same irrespective of whether or not $\lambda^0$ is known. In addition, they proved that the estimate of $\lambda^0$ is $n$-consistent although its limit distribution remains an open problem.

However, if the ultimate interest is only the estimation of the memory parameter $\alpha$, one possible criticism of the parametric approach is that an incorrect specification of the model leads to inconsistent estimates of $\alpha$. One source of misspecification is the choice of a wrong value of $\lambda^0$. If that were the case, Whittle estimates of $\alpha$ would be inconsistent, and would possibly estimate the value *zero*. The latter might happen even if a semiparametric approach were adopted; see Section 3. Thus, we might conclude that the data is short-memory- instead of long-memory-dependent, which could have some adverse effects on the statistical inference of relevant statistics such as the serial covariances; see [21] or [33].

The main objectives of this paper are twofold: first, under mild conditions, to provide a consistent estimator of $\lambda^0$ and characterize its limit distribution. In particular we show that the limiting distribution of the estimator of $\lambda^0$ depends on whether $\lambda^0 \in (0,\pi)$ or $\lambda^0 = 0$ or $\pi$. The second objective is to investigate the consequences that the lack of knowledge of $\lambda^0$ might have on the estimation of $\alpha$.

Some earlier related work has been completed by Müller and Prewitt [27] and Yajima [36]. In the former the authors estimate the peak of the spectral density $f(\lambda)$ in a model, like that in (1.2), whose spectral density function is continuous in $[0,\pi]$. Looking at $\arg\sup_\lambda \hat{f}(\lambda)$, where $\hat{f}(\lambda)$ is a smoothed nonparametric estimate of $f(\lambda)$, they show its consistency and the limit distribution to be a normal random variable when $\lambda^0 \in (0,\pi)$. Yajima [36] considers the estimation of $\lambda^0$ in a model with spectral density function possessing a pole at $\lambda^0$. Based on the maximum of the periodogram of the data, he gives consistency and an upper rate of convergence for the estimate of $\lambda^0$. Unfortunately, the limit distribution, which is required for statistical inference, was not provided. In addition, his results rely on the assumption that the data is Gaussian, which is not required in the present paper. Finally, it should be mentioned that Giraitis and Leipus [13] prove the consistency of $\lambda^0$ in a model like (1.6).

The paper is organized as follows. In Section 2 we describe a semiparametric estimator of $\alpha$ when $\lambda^0$ and the estimator $\hat{\lambda}^0$ of $\lambda^0$ are known. In Section 3 we discuss the statistical properties of $\hat{\lambda}^0$ and we show that the asymptotic properties of a two-step estimator of $\alpha$ remain the same irrespective of whether $\lambda^0$ is known or estimated. The finite sample behavior of the estimators of $\lambda^0$ and $\alpha$ is analyzed in Section 4 through a Monte Carlo study. Section 5 provides the proofs of the results given in Section 3, which apply a series of lemmas given in Section 6. Finally, Section 7 contains a summary.

**2. Regularity conditions and the estimators of the pole and memory parameter.** Let $\{x_t\}$ be a covariance stationary linear process observed at times $t = 1, 2, \ldots, n$, with spectral density $f(\lambda)$ satisfying (1.4). When $\lambda^0$ is known, under the semiparametric specification (1.4) several estimators of the memory parameter $\alpha$ have been proposed and examined. In this paper, to estimate $\alpha$ we shall use a modification of the log-periodogram estimator (see [31]), which we now describe. Consider the average periodogram spectral density estimator of $f(\lambda)$,

$$\ddot{f}_\ell = \ddot{f}(\lambda_\ell) = \frac{1}{2k_1 + 1} \sum_{|j| \leq k_1} I_{\ell+j}, \tag{2.1}$$

where $I_\ell = I(\lambda_\ell)$ denotes the periodogram of $x_t$, that is,

$$I_\ell = \left|(2\pi n)^{-1/2} \sum_{t=1}^n x_t e^{it\lambda_\ell}\right|^2, \qquad \ell = 1, \ldots, [n/2], \tag{2.2}$$

and $I_0 = 0$, where $\lambda_\ell = (2\pi\ell)/n$ for $\ell = 1, 2, \ldots, [n/2]$, $[z]$ denotes the integer part of $z$ and $k_1 = k_1(n)$ is a positive number which increases slowly with $n$, that is, $k_1^{-1} + n^{-1}k_1 \to 0$. (Observe that the definition of $I_\ell$ entails sample-mean correction.)

Let $\psi(u)$ be a weight function in $(0,1)$, and write

$$(2.3) \qquad \hat{\alpha}(\lambda_q) = \frac{1}{2\bar{h}_\psi k} \sum_{p=1}^{k} \psi_p (\log \hat{f}_{q+p} + \log \hat{f}_{q-p}),$$

where $\psi_p = \psi(p/k)$, $\bar{h}_\psi = -k^{-1} \sum_{p=1}^{k} \psi_p \log(p/k)$, $\hat{f}_\ell = \max(\ddot{f}_\ell, n^{-1})$ and $k = k(n)$, a positive number which increases slowly with $n$, that is, $k^{-1} + n^{-1}k \to 0$.

DEFINITION 2.1. If $\lambda^0$ is known, we define the estimator of the memory parameter $\alpha$ as $\hat{\alpha}(\lambda_s)$, where $\lambda_s$ is the closest Fourier frequency $\lambda_q$ to $\lambda^0$.

REMARK 2.1. The motivation to use $\hat{f}_\ell$ instead of $\ddot{f}_\ell$ in (2.3) is due to the singular behavior of $\log x$ at $x = 0$. Specifically for the proof of tightness, that is, Proposition 5.4 in Section 5, we have not been able to bound some probabilities or moments for all $n \geq n_0$ as required. This problem, of course, does not appear as $n \to \infty$ as can be observed from Propositions 5.1–5.3, nor if our goal were to examine the behavior of $\hat{\alpha}(\lambda_s)$. We do not believe that this adjustment is needed in practice and have made it here only because we cannot establish Theorem 3.2 (cf. Proposition 5.4) without it, unless some additional stronger conditions were introduced, for instance, the normality of the data.

We now define our estimator of $\lambda^0$ as $\hat{\lambda}^0 = \lambda_{\hat{q}} = (2\pi\hat{q})/n$, where

$$(2.4) \qquad \hat{q} = \underset{q=0,\ldots,[n/2]}{\arg\max} \ \hat{\alpha}(\lambda_q).$$

Note that periodicity and symmetry around zero imply that it suffices to search for the maximum in (2.4) at frequencies $\lambda_q$, with $q = 0, \ldots, [n/2]$. From (2.4) we could define $\hat{\alpha}(\hat{\lambda}^0)$ as an estimator of $\alpha$, that is, (2.3) evaluated at $\hat{\lambda}^0$. However (see Section 3), since $\hat{\alpha}(\hat{\lambda}^0)$ does not have optimal properties, we will describe a two-step estimator, denoted $\check{\alpha}(\check{\lambda}^0)$, which overcomes all the adverse properties of $\hat{\alpha}(\hat{\lambda}^0)$; see Theorem 3.4.

The motivation for the estimator in (2.4) is as follows. From the proof of Theorem 3.4(a) below, it is easily shown that $\hat{\alpha}(\lambda_s)$ is a consistent estimator of $\alpha$. On the other hand, if $\lambda_q$ is in any open set outside $\lambda^0$, that is, $|\lambda_q - \lambda_s| > \delta$, for any arbitrarily small $\delta > 0$, Condition C.1 below implies that

$$f(\lambda_q) = |\lambda_q - \lambda^0|^{-\alpha} g(\lambda_q) \sim C g(\lambda_q).$$

That is, at the frequency $\lambda_q$ the spectral density function behaves as if $\alpha$ were equal to zero. So, from the proof of Theorem 3.4(a) we should expect that $\hat{\alpha}(\lambda_q) \xrightarrow{P} 0$, implying that $\Pr\{|\hat{\lambda}^0 - \lambda_s| < \delta\} \to 1$. That is, the estimator given in (2.4) is consistent. These heuristics will be formalized in Theorem 3.1 below. We finish this section by introducing the following regularity conditions and their discussion.

CONDITION C.1. There exists $\alpha \in (0,1)$ such that
$$f(\lambda) = \begin{cases} |\lambda - \lambda^0|^{-\alpha} g(\lambda), & \text{if } 0 \leq \lambda \leq \pi, \\ |\lambda + \lambda^0|^{-\alpha} g(\lambda), & \text{if } -\pi \leq \lambda \leq 0, \end{cases}$$
where $g(\lambda)$ is a bounded symmetric and bounded away from zero function with two continuous derivatives for $0 < \lambda < \pi$.

CONDITION C.2. $\{x_t\}$ is a covariance stationary linear process,
$$x_t = \sum_{j=0}^{\infty} \beta_j \varepsilon_{t-j}, \qquad \sum_{j=0}^{\infty} \beta_j^2 < \infty,$$
where $\{\varepsilon_t\}$ is a zero-mean i.i.d. sequence with $E(\varepsilon_t^2) = 1$ and $E|\varepsilon_t|^\ell = \mu_\ell < \infty$ for $\ell = 3, \ldots, 2\tau$ and some $\tau \geq 4$.

CONDITION C.3. As $\lambda \to \lambda$, the function $\beta(\lambda) = \sum_{j=0}^{\infty} \beta_j e^{ij\lambda}$ satisfies
$$|\partial \beta(\lambda)/\partial \lambda| = O(|\lambda - \lambda^0|^{-1} |\beta(\lambda)|).$$

CONDITION C.4. $k^{1+\iota} k_1^{-2} + k^{-2} k_1^3 \log k + n k_1^{-(\tau^2+2)/2\tau} \to 0$, for some $\iota > 0$ as $n \to \infty$, with $k \leq c n^{4/5}$, $0 < c < \infty$ and where $\tau$ is as in Condition C.2.

CONDITION C.5. The function $\psi(x)$ is twice continuously differentiable with second derivative that is Lipschitz of order at least $\frac{1}{2}$ in its support $(0,1)$ and satisfies $\int_0^1 \psi(x)\,dx = 0$, $0 < h_\psi = -\int_0^1 \psi(x)(\log x)\,dx < \infty$, $0 < \bar{\psi}'' = \int_0^1 \psi''(x)(\log x)\,dx$, where $\psi''(x) = \frac{d^2}{dx^2}\psi(x)$, and $|x^{-2}\psi(x)| + |(1-x)^{-1}\psi(x)| \leq D < \infty$.

We now discuss Conditions C.1–C.5. Condition C.1 is much the same as that employed by Robinson [31, 32]. Indeed, Condition C.1 implies that as $\lambda \to \lambda^0$,
$$f(\lambda) = C|\lambda - \lambda^0|^{-\alpha}(1 + C^{-1} g'(\lambda^0)(\lambda - \lambda^0) + O(|\lambda - \lambda^0|^2))$$
by Taylor expansion of $g(\lambda)$ around $\lambda^0$ and where $C = g(\lambda^0)$. Observe that $g'(\lambda^0) = 0$ when $\lambda^0 = \{0, \pi\}$ by symmetry of $f(\lambda)$, obtaining then the corresponding condition used in [31, 32]. However, we prefer to state the condition in its present form since, in Theorems 3.1 and 3.2 below, some regularity conditions on $f(\lambda)$ are needed outside any open set containing $\pm\lambda^0$. Examples of processes whose spectral density function satisfies Condition C.1 are the FARIMA$(p, \alpha/2, q)$ and the GARMA model given in (1.6). Finally, the last part of Condition C.1 is quite standard in the spectral density estimation literature. Condition C.2 is needed for the proof of tightness (see the proofs of Theorems 3.2 and 3.1). It is also required to show the uniform convergence of $\hat{f}$, although for the latter property, at the expense of stronger conditions on the rate of convergence of $k_1^{-1}$ to zero, fewer moments of $\varepsilon_t$ can be assumed. Obviously Condition C.2 is satisfied if $\varepsilon_t$ is Gaussian. Condition C.3 is the same as Robinson's [32]. Condition C.4 controls the rate of increase of $k$ and $k_1$. For instance, denoting $k = n^{\gamma_2}$ and $k_1 = n^{\gamma_1}$, in the Gaussian case, we can take $0 < \gamma_1 < 8/15$ and $3\gamma_1/2 < \gamma_2 \leq \min\{2\gamma_1/(1+\iota), 4/5\}$, whereas for $\tau = 4$ the bounds are $4/9 < \gamma_1 < 8/15$ and $3\gamma_1/2 < \gamma_2 \leq 4/5$. Finally, Condition C.5 characterizes the type of weight in (2.3). An example of $\psi(u)$ satisfying Condition C.5 is $\psi(u) = -u^2 + 35u^{2.5}/6 - 29u^3/6 + 2u^3 \log u$.

It is worth mentioning that the quadratic behavior of the weight $\psi(u)$, as $u \to 0$, guarantees that the first moment of $\xi(v)$ (see Theorem 3.2 below) has a parabolic structure so that the maximum of $\xi(v)$ is easily obtainable. Obviously, other different types of weights can be used which would not prevent the consistency of the estimator of $\lambda^0$. However, for weights not having a quadratic behavior, the asymptotic distribution of the estimate of the pole is not guaranteed to be normally distributed. We will return to this condition after Theorem 3.2.

**3. Statistical properties of the estimators of the pole and memory parameter.** In this section we prove a functional limit theorem for a process operating on increments of $\hat{\alpha}(\lambda_q)$ near $\lambda^0$, which together with the continuous mapping theorem will allow obtaining the asymptotic distribution of $\hat{\lambda}^0$. A similar approach was used by Eddy [10] to estimate the mode of a probability density function and by Müller [26] for the estimation of the break point in a regression model. Apart from providing the consistency and rate of convergence of $\hat{\lambda}^0$ to $\lambda^0$, the limit distribution will guarantee that asymptotic valid inferences around the true value of $\lambda^0$ may be implemented.

The strategy of the proof to obtain the asymptotic distribution of $\hat{\lambda}^0$ consists of three steps; see [34], Chapter 3. Step 1 establishes the consistency of $\hat{\lambda}^0$ to $\lambda^0$. Step 2 establishes the rate of convergence of $\hat{\lambda}^0$ to $\lambda^0$, and Step 3 shows that suitably rescaled versions of $\hat{\alpha}(\lambda_q)$ converge weakly to a limit, denoted $\xi(v)$ in Theorem 3.2, in the space $\mathbb{D}[-M, M]$ for each finite $0 < M < \infty$. Note that convergence in $\mathbb{D}[-M, M]$ for each finite $0 < M < \infty$ is to be meant convergence in $\mathbb{D}(-\infty, \infty)$. See Pollard [28]. From here, the continuous mapping theorem will conclude that $\hat{\lambda}^0$, after normalization, will converge in distribution to the $\arg\max_v \xi(v)$.

The next theorem gives the consistency and rate of convergence of $\hat{\lambda}^0$ to $\lambda^0$, that is, Steps 1 and 2. Theorem 3.2 justifies Step 3, whereas Corollary 3.3 examines the asymptotic distribution of $\hat{\lambda}^0$.

THEOREM 3.1. *Assuming Conditions* C.1–C.5, $|\hat{\lambda}^0 - \lambda^0| = O_p(k^{1/2} n^{-1})$.

We see that the rate of convergence of $\hat{\lambda}^0$ to $\lambda^0$ is slower than the parametric rate $n^{-1}$ obtained by Giraitis, Hidalgo and Robinson [12]. This appears to be reasonable due to the local behavior of our statistics. The same phenomenon occurs in other related, although different, problems involving nonparametric statistics; see, for instance, [10, 26] or [27].

Under Conditions C.2 and C.4, for $\tau = 4$, $(\hat{\lambda}^0 - \lambda^0) = O_p(n^{\delta - 2/3})$ for any arbitrarily small $\delta > 0$. However, a closer examination of these conditions and the proof of Lemma 6.3 indicate that the rate depends on the number of finite moments of the sequence $\varepsilon_t$ in Condition C.2. In general, with $\tau \geq 4$, $(\hat{\lambda}^0 - \lambda^0) = O_p(n^{\delta - (2\tau^2 - 3\tau + 4)/2(\tau^2 + 2)})$. So, the greater the number of finite moments allowed for $\varepsilon_t$, the faster the rate of convergence of $\hat{\lambda}^0$ to $\lambda^0$. In the extreme case where all the moments exist, the rate of convergence of $\hat{\lambda}^0$ becomes $n^{\delta - 1}$. This rate was obtained by Yajima [36] in the Gaussian case and is arbitrarily close to $n^{-1}$ obtained in [12].

So Theorem 3.1 indicates that $\hat{\lambda}^0 = \lambda^0 + n^{-1}(2\pi[k^{1/2}v])$ for some $|v| \leq M < \infty$. To examine the asymptotic distribution of $\hat{\lambda}^0$, let us introduce the notation

$$(3.1) \qquad \hat{\xi}_n(v) = k(\hat{\alpha}(\lambda_{s+[k^{1/2}v]}) - \hat{\alpha}(\lambda_s)).$$

$\hat{\xi}_n(v)$ is a random step function which is constant in the intervals $[i/k^{1/2}, (i+1)/k^{1/2})$, $|i| \leq M$, so that $\hat{\xi}_n(v)$ is a random element in the Skorohod space $\mathbb{D}[-M, M]$ for arbitrary $0 < M < \infty$.

We now establish our main result, that is, the aforementioned Step 3.

THEOREM 3.2. *Let $|v| \leq M$ for any arbitrary $M \in (0, \infty)$. Assuming Conditions* C.1–C.5,

$$\hat{\xi}_n(v) \overset{weakly}{\Longrightarrow} \xi(v) \qquad \text{in the space } \mathbb{D}[-M, M],$$

*where $\xi(v)$ is a continuous Gaussian process such that*

$$E(\xi(v)) = -h_\psi^{-1} \bar{\psi}'' v^2 \alpha / 2 \quad \text{and} \quad \mathrm{Cov}(\xi(v_1), \xi(v_2)) = h_\psi^{-2} \varsigma v_1 v_2,$$

*where $\varsigma = \int_0^1 \psi'(u)^2 \, du < \infty$ with $\psi'(x) = \frac{d}{dx} \psi(x)$.*

The immediate consequence of Theorem 3.2 is that $\arg\max_v \xi(v)$ is a normal random variable. Indeed, because Theorem 3.2 implies that the limiting process $\xi(v)$ is Gaussian, it can be written as

$$\xi(v) = -h_\psi^{-1}\bar{\psi}'' v^2 \alpha/2 + h_\psi^{-1} \varsigma^{1/2} v X,$$

where $X = N(0,1)$. But $\xi(v)$ is a random parabola with fixed second derivatives and a unique maximum at

$$v^* = (\bar{\psi}''\alpha)^{-1} \varsigma^{1/2} X,$$

since by Condition C.5, $0 < h_\psi, 0 < \bar{\psi}''$, so that $\partial^2 \xi(v)/\partial v^2 = -h_\psi^{-1}\bar{\psi}''\alpha < 0$. From here we can observe the (possible) consequences of using a weight function $\psi(u)$ which does not have a parabolic structure at $u=0$. The main implication is that if the latter were the case, $E(\xi(v))$ would not necessarily be a parabola as in Theorem 3.2. For example, it may be that $E(\xi(v)) = C|v|$, in which case not only can the $\arg\max_v \xi(v)$ be difficult to obtain, but more importantly it would no longer be a normal random variable. So, in view of the asymptotic normality achieved with a weight $\psi(u)$ satisfying Condition C.5, it appears desirable to employ it. Similar issues occur when estimating the date of a break in a regression model; see, for example, [26].

Now we turn our attention to the asymptotic properties of $\hat{\lambda}^0$. Note that Theorem 3.1 indicates that

$$(3.2) \qquad \hat{\lambda}^0 = \lambda_s + \frac{2\pi k^{1/2}}{n}\hat{v}_n = \lambda^0 + \frac{2\pi k^{1/2}}{n}\hat{v}_n + O\left(\frac{1}{n}\right),$$

where $\hat{v}_n = \arg\max_v \hat{\xi}_n(v)$. Then we have the following:

COROLLARY 3.3. *Denote $\Psi = \varsigma(\bar{\psi}''\alpha)^{-2}$. Assuming Conditions C.1–C.5, as $n \to \infty$:*

(a) *If $\lambda^0 \in (0,\pi)$, then $(2\pi k^{1/2})^{-1} n(\hat{\lambda}^0 - \lambda^0) \xrightarrow{d} Z_{\lambda^0} \equiv Y = N(0, \Psi)$.*

(b) *If $\lambda^0 = 0$, then $(2\pi k^{1/2})^{-1} n\hat{\lambda}^0 \xrightarrow{d} Z_0 = Y\mathcal{I}(Y \geq 0)$, where $\mathcal{I}(A)$ denotes the indicator function of the set $A$.*

(c) *If $\lambda^0 = \pi$, then $(2\pi k^{1/2})^{-1} n(\hat{\lambda}^0 - \pi) \xrightarrow{d} Z_\pi = Y\mathcal{I}(Y \leq 0)$.*

We now comment on the results of Corollary 3.3. First, we now see the necessity of Theorem 3.1, as it will give us the normalization needed to achieve a "proper" asymptotic distribution. Next, we observe that the limiting distribution of $\hat{\lambda}^0$ depends on whether $\lambda^0$ is $\{0,\pi\}$ or $\lambda^0 \in (0,\pi)$. The intuition about the limiting distribution of $\hat{\lambda}^0$ in cases (b) and (c) is as follows. As the maximization of $\hat{\alpha}(\lambda_q)$ in (2.4) is restricted to the interval $0 \leq \lambda_q \leq \pi$, for $\lambda^0 = 0$, it implies that $(\hat{\lambda}^0 - \lambda^0) = \hat{\lambda}^0 \geq 0$ so that $Z_0$ cannot take negative values. Similarly, $\lambda^0 = \pi$ implies that $\hat{\lambda}^0 - \pi \leq 0$ and $Z_\pi$ cannot take positive values. So, the estimation of $\lambda^0$ falls into the category of a constrained optimization problem or inequality constraint estimation. Indeed, when $\lambda^0$ is an interior point of the set $[0,\pi]$ and due to the consistency of $\hat{\lambda}^0$, we can expect that the constrained estimator, $\hat{\lambda}^0 = \lambda_{\hat{q}}$, coincides with the unconstrained estimator $\tilde{\lambda}_q = \lambda_{\tilde{q}} = (2\pi\tilde{q})/n$, where

$$(3.3) \qquad \tilde{q} = \underset{q \in \{0, \pm 1, \pm 2, \dots\}}{\arg\max} \hat{\alpha}(\lambda_q),$$

whereas if $\lambda^0 = 0$, $\hat{q} = \tilde{q}\mathcal{I}(\tilde{q} \geq 0)$. Similar arguments apply when $\lambda^0 = \pi$.

Once we have examined the properties of $\hat{\lambda}^0$, we next examine the estimation of $\alpha$. By Theorems 3.1 and 3.2 and the functional mapping theorem, it is easily shown that $\hat{\alpha}(\hat{\lambda}^0) - \hat{\alpha}(\lambda^0) = o_p(k^{-1/2})$. So, $k^{1/2}(\hat{\alpha}(\hat{\lambda}^0) - \alpha)$ and $k^{1/2}(\hat{\alpha}(\lambda^0) - \alpha)$ have the same asymptotic distribution. However, the faster the convergence of $\hat{\lambda}^0$ to $\lambda^0$, the slower the rate of convergence of $\hat{\alpha}(\hat{\lambda}^0)$ to $\alpha$, and hence it becomes slower than the rate obtained when $\lambda^0$ is known. The same phenomenon happens to hold in [26]. Hence, to circumvent this drawback, as in [26], we propose a *two-step* procedure to estimate $\alpha$. To this end, we shall use as an estimator of $\alpha$ that given in (2.3) where $\lambda_q$ is replaced by $\breve{\lambda}^0 = (2\pi\breve{q})/n$ such that $|\breve{\lambda}^0 - \lambda^0| = O_p(k^{1/2}/n)$, and $k$ and $k_1$ are replaced by $m$ and $m_1$, respectively, satisfying:

CONDITION C.6. $m^{-1} + mm_1^{-2} + m_1^5 m^{-3} \log m_1 + k/m \to 0$ with $m = cn^{4/5}$, $0 < c < \infty$.

In addition, to be a bit more general regarding our choice of the weight function $\psi(u)$, we allow for the weighted function, say $w(u)$, to satisfy:

CONDITION C.7. $\int_0^1 w(u)\,du = 0$, $0 < h_w = -\int_0^1 w(u)(\log u)\,du < \infty$, $w(u) \sim cu^\zeta$ as $u \to 0+$ for some $1/3 \leq \zeta \leq 1$ and for all $0 < u_1 < u_2 < 1$,

$$|w(u_2) - w(u_1)| \leq D|u_2 - u_1|^\zeta, \qquad 0 < D < \infty.$$

So, our *two-step* estimator of $\alpha$ is defined as

$$(3.4) \qquad \breve{\alpha}(\breve{\lambda}^0) = \frac{1}{2\bar{h}_w m} \sum_{p=1}^m w_p (\log \hat{f}_{\breve{q}+p} + \log \hat{f}_{\breve{q}-p}),$$

where $\hat{f}(\lambda) = \max\{\ddot{f}(\lambda), n^{-1}\}$ and $\ddot{f}(\lambda)$ is as in (2.1) but with the smoothing parameter $k_1$ there being replaced by $m_1$, $\bar{h}_w = -m^{-1}\sum_{p=1}^m w_p \log(p/m)$ and $w_p = w(p/m)$.

We now comment on $\breve{\alpha}(\breve{\lambda}^0)$ compared to $\ddot{\alpha}(\breve{\lambda}^0) = (\bar{h}_w m)^{-1} \sum_{p=1}^m w_p \log \hat{f}_{\breve{q}+p}$. Observe that the former is a "symmetrized" version of the latter $\ddot{\alpha}(\breve{\lambda}^0)$. Assume for simplicity that $\lambda^0$ is known. As in other semiparametric estimators, for example, [31], one source of the bias of $m^{1/2}(\ddot{\alpha}(\lambda^0) - \alpha)$ comes from the replacement of $f(\lambda)$ by $g(\lambda^0)|\lambda - \lambda^0|^{-\alpha}$, which in our case, that is, if $\lambda^0 \neq \{0, \pi\}$, will be proportional to

$$m^{-1/2} \sum_{p=1}^{m} w_p(\lambda_p + O(\lambda_p^2)) = O(n^{-1}m^{3/2}).$$

The main reason for this behavior is that when $\lambda^0 = \{0, \pi\}$, by symmetry we have $g'(0) = g'(\pi) = 0$, whereas for $\lambda \neq 0$ or $\pi$, $g'(\lambda)$ may not be zero so that $g^{-1}(\lambda^0)|\lambda - \lambda^0|^\alpha f(\lambda) = 1 + g^{-1}(\lambda^0)g'(\lambda^0)(\lambda - \lambda^0) + O(|\lambda - \lambda^0|^2)$ by a Taylor expansion of $g(\lambda)$ around $\lambda^0$. Recall the comments made on Condition C.1. However, when the estimator $\check{\alpha}(\lambda^0)$ in (3.4) is employed, the contribution of the above approximation (Taylor expansion) to the bias of $m^{1/2}(\check{\alpha}(\lambda^0) - \alpha)$ is proportional to

$$m^{-1/2} \sum_{p=1}^{m} w_p(-\lambda_p + O(\lambda_p^2)) + m^{-1/2} \sum_{p=1}^{m} w_p(\lambda_p + O(\lambda_p^2)) = O(m^{5/2}n^{-2}).$$

Note that the latter holds true also for $\lambda^0 = \{0, \pi\}$. So, the "symmetrized" estimator $\check{\alpha}(\lambda^0)$ would have a smaller bias order and thus it would have a faster rate of convergence to $\alpha$ than $\ddot{\alpha}(\lambda^0)$.

THEOREM 3.4. *Denote* $\Phi^2 = 2^{-1} \int_0^1 w^2(x)\,dx$ *and* $B = (\partial^2/\partial\lambda^2 \log g(\lambda^0)) \times \int_0^1 u^2 w(u)\,du$. *Let* $\check{\lambda}^0$ *be an estimator of* $\lambda^0$ *such that* $|\check{\lambda}^0 - \lambda^0| = O_p(k^{1/2}/n)$. *Assuming Conditions* C.1–C.4 *with* $k_1 = n^{\gamma_1}$ *and* $k = n^{\gamma_2}$, *where* $2\tau/(\tau^2 + 2) < \gamma_1 < 8/15$, $3\gamma_1/2 < \gamma_2 < \min\{\frac{2\gamma_1}{1+\iota}, \frac{4}{5}\}$, $\tau$ *as in Condition* C.2, *and Conditions* C.6 *and* C.7, *then*

(a) $(2m)^{1/2}(\check{\alpha}(\lambda_s) - \alpha) \xrightarrow{d} N(4\pi^2 c^{5/2}B/(2^{1/2}h_w), \Phi^2/h_w^2)$,
(b) $(2m)^{1/2}(\check{\alpha}(\check{\lambda}^0) - \alpha) \xrightarrow{d} N(4\pi^2 c^{5/2}B/(2^{1/2}h_w), \Phi^2/h_w^2)$.

REMARK 3.1. It is worth mentioning that the results of Theorem 3.4 hold true if the weight $\psi(u)$ employed to estimate $\lambda^0$ is used in $\check{\alpha}(\lambda^0)$. However, this weight will not guarantee an asymptotic variance smaller than 1, as is the case with the weight used in the Monte Carlo experiment. In fact, for the weight used in the Monte Carlo experiment, $h_w^{-2}\Phi^2 \sim 0.70$, which is smaller than the corresponding asymptotic variance of other estimators of $\alpha$ suggested in the literature. Finally, the theorem indicates that although any preliminary estimator of $\lambda^0$ which satisfies $|\check{\lambda}^0 - \lambda^0| = O_p(k^{1/2}/n)$ is adequate for the results to follow, in practice it appears that one may use that given in (2.4) for computational simplicity.

Theorem 3.4 provides a consistent estimator of the asymptotic variance of $\hat{\lambda}^0$, that is, $\Psi$ in Corollary 3.3, by replacing $\alpha$ by $\check{\alpha}(\check{\lambda}^0)$. But more importantly, it indicates that the *two-step* estimator $\check{\alpha}(\check{\lambda}^0)$, apart from having the same asymptotic distribution as $\check{\alpha}(\lambda^0)$, achieves the optimal semiparametric rate of convergence obtained by Giraitis, Robinson and Samarov [14] when $\lambda^0 = 0$. So, asymptotically, there is no loss by using $\hat{\lambda}^0$ instead of $\lambda^0$. However, to achieve the latter, as in other nonparametric estimates, $\check{\alpha}(\hat{\lambda}^0)$ will have a bias term of the same order of magnitude as the standard deviation.

**4. Finite sample behavior.** In this section we study via Monte Carlo analysis the finite sample performance of the estimators $\hat{\lambda}^0$ and $\check{\alpha}(\hat{\lambda}^0)$. The models employed throughout the simulations are

$$(4.1) \qquad (1-L)^{\alpha/2} x_t = \varepsilon_t, \qquad t = 0, \pm 1, \ldots,$$

$$(4.2) \qquad (1 - 2\cos(\pi/2)L + L^2)^{\alpha/2} x_t = \varepsilon_t, \qquad t = 0, \pm 1, \ldots,$$

where $\{\varepsilon_t\}$ is a zero-mean sequence of i.i.d. Gaussian random variables. Model (4.1) generates a pole at $\lambda^0 = 0$, whereas model (4.2) does so at $\lambda^0 = \pi/2$. We have chosen $\alpha = 0.2, 0.4, 0.6$ and $0.8$. The autocorrelation functions of (4.1) and (4.2) are given by

$$\rho_j = \frac{j - 1 + \alpha/2}{j - \alpha/2} \rho_{j-1}, \qquad j = 1, 2, \ldots,$$

and

$$\rho_{2j} = \frac{1 - j - \alpha/2}{j - \alpha/2} \rho_{2(j-1)}, \qquad \rho_{2j-1} = 0, \qquad j = 1, 2, \ldots,$$

respectively; see, for example, [3]. For each combination of $\alpha$ and $\lambda^0$, 2500 replications of series of lengths $n = 256$ and $1024$ were generated by the method of Davies and Harte [8].

Also, we have compared the performance of $\hat{\lambda}^0$ and $\check{\alpha}(\hat{\lambda}^0)$ with the corresponding estimators obtained using the log-periodogram estimator of [31] popular among practitioners. That is, consider $\tilde{\lambda}^0 = \lambda_{\tilde{q}} = (2\pi\tilde{q})/n$, where

$$(4.3) \qquad \tilde{q} = \underset{q=0,\ldots,[n/2]}{\arg\max}\ \hat{\alpha}_{\mathrm{LOG}}(\lambda_q),$$

$$(4.4) \qquad \hat{\alpha}_{\mathrm{LOG}}(\lambda_q) = -\left(2 \sum_{j=1}^{k} \phi_j \log j\right)^{-1} \sum_{j=1}^{k} \phi_j (\log I_{j+q} + \log I_{q-j}),$$

with $\phi_j = \log j - k^{-1} \sum_{\ell=1}^{k} \log \ell$. Moreover, we have examined the behavior of the estimator $\hat{\alpha}_{\mathrm{LOG}}(\tilde{\lambda}^0)$ of $\alpha$, where $k = m$ in (4.4). For the estimation of $\lambda^0$, the chosen bandwidth parameters were, for $n = 256$ and $1024$, $k = 14$ and $24$, respectively, and $k_1 = k^{0.6} \log \log 2k$, whereas for the *two-step* estimators $\check{\alpha}(\hat{\lambda}^0)$ and $\hat{\alpha}_{\mathrm{LOG}}(\tilde{\lambda}^0)$ of $\alpha$, we have chosen $m = n/4$ and $m_1 = m^{0.6} \log \log 2m$. The weight functions used were $\psi(u) = -u^2 + 35u^{2.5}/6 - 29u^3/6 + 2u^3 \log u$ and $w(u) = u^{1/3} - 9u^{1/2}/8$, respectively.

SEMIPARAMETRIC ESTIMATION OF THE POLE 13

Table 1 illustrates the bias and standard deviation of the estimators $\hat{\lambda}^0$ given in (2.4) and $\tilde{\lambda}^0$ in (4.3). More specifically, since $\hat{\lambda}^0 = (2\pi\hat{q})/n$ and $\tilde{\lambda}^0 = (2\pi\tilde{q})/n$, we have reported the bias and standard deviation of $\hat{q}$ and $\tilde{q}$. Table 2 summarizes the bias, standard deviation and mean square error of $\breve{\alpha}(\hat{\lambda}^0)$ and $\breve{\alpha}(\lambda^0)$. The motivation to include $\breve{\alpha}(\lambda^0)$ is to investigate the relative loss we incur by lack of knowledge of $\lambda^0$ in small samples. Recall that Theorem 3.4 indicates that asymptotically there is no loss. Moreover, Table 2 illustrates the finite sample performance of the corresponding estimators of $\alpha$ obtained using the log-periodogram estimator in (4.4), that is, $\hat{\alpha}_{\text{LOG}}(\lambda^0)$ and $\hat{\alpha}_{\text{LOG}}(\tilde{\lambda}^0)$.

Inspection of Table 1 indicates better performance of $\hat{\lambda}^0$ than $\tilde{\lambda}^0$ across different models and sample sizes, especially for $\alpha > 0.2$. For example, when $\alpha = 0.8$, the finite sample performance of $\hat{\lambda}^0$ is clearly superior to that of $\tilde{\lambda}^0$, this superiority being greater with the sample size. With regard to the estimators of the memory parameter $\alpha$, we observe that the proposed two-step estimator $\breve{\alpha}(\hat{\lambda}^0)$ outperforms $\hat{\alpha}_{\text{LOG}}(\tilde{\lambda}^0)$ and has better finite sample properties for all $\alpha$ and $\lambda^0$. In some cases, the performance of $\hat{\alpha}_{\text{LOG}}(\tilde{\lambda}^0)$ is very poor compared to that of $\breve{\alpha}(\hat{\lambda}^0)$, especially for large values of $\alpha$. Finally, when comparing their performances with the estimators obtained when the location of the pole $\lambda^0$ is known, we observe that the relative loss of efficiency of $\breve{\alpha}(\hat{\lambda}^0)$ is smaller than that of $\hat{\alpha}_{\text{LOG}}(\tilde{\lambda}^0)$. Moreover, as Theorem 3.4 indicates, it appears that knowledge of $\lambda^0$ is not relevant to estimate $\alpha$ when $\breve{\alpha}(\hat{\lambda}^0)$ is used, although it seems not to be the case when the log-periodogram is employed. Altogether, we can conclude that $\hat{\lambda}^0$ and $\breve{\alpha}(\hat{\lambda}^0)$ enjoy better finite sample properties than the corresponding ones based on $\tilde{\lambda}^0$ and $\hat{\alpha}_{\text{LOG}}(\tilde{\lambda}^0)$.

TABLE 1
*Bias and standard deviation of $\hat{q}$ and $\tilde{q}$*

|  |  |  | $\alpha$ | | | | | | | |
|---|---|---|---|---|---|---|---|---|---|---|
|  |  |  | 0.2 | | 0.4 | | 0.6 | | 0.8 | |
| $\lambda^0 = 0$ | $n$ | 256 | 9.35 | (8.33) | 6.38 | (6.96) | 4.24 | (5.39) | 2.80 | (4.04) |
|  |  |  | 9.26 | (7.88) | 7.32 | (6.85) | 5.94 | (6.01) | 4.85 | (5.25) |
|  |  | 1024 | 15.40 | (15.50) | 8.43 | (10.74) | 4.81 | (7.64) | 2.62 | (5.76) |
|  |  |  | 22.91 | (25.31) | 15.55 | (20.89) | 9.60 | (14.02) | 6.73 | (9.96) |
| $\lambda^0 = \frac{\pi}{2}$ | $n$ | 256 | 0.003 | (7.64) | $-0.084$ | (5.33) | $-0.091$ | (2.96) | $-0.054$ | (1.56) |
|  |  |  | 0.209 | (9.59) | 0.270 | (9.21) | 0.272 | (8.66) | 0.320 | (7.28) |
|  |  | 1024 | 0.051 | (11.87) | 0.117 | (4.77) | 0.063 | (1.89) | 0.216 | (1.13) |
|  |  |  | 0.435 | (27.89) | 0.144 | (25.77) | $-0.213$ | (21.30) | $-0.060$ | (13.52) |

The first row in each cell corresponds to $\hat{q}$, whereas the second row is that of $\tilde{q}$.

Table 2
*Bias, standard deviation and MSE of the long-memory parameter estimators*

| | | \multicolumn{12}{c}{$\alpha$} |
| | | 0.2 | | | 0.4 | | | 0.6 | | | 0.8 | | |
| $\lambda^0$ | $n$ | BIAS | S.D. | M.S.E. | BIAS | S.D. | M.S.E. | BIAS | S.D. | M.S.E. | BIAS | S.D. | M.S.E. |
|---|---|---|---|---|---|---|---|---|---|---|---|---|---|
| 0 | 256 | −0.020 | 0.064 | 0.004 | −0.022 | 0.067 | 0.005 | −0.017 | 0.071 | 0.005 | −0.006 | 0.072 | 0.005 |
| | | −0.019 | 0.057 | 0.004 | −0.030 | 0.065 | 0.005 | −0.024 | 0.074 | 0.006 | −0.006 | 0.075 | 0.006 |
| | | −0.001 | 0.089 | 0.008 | −0.003 | 0.089 | 0.008 | −0.003 | 0.089 | 0.008 | −0.007 | 0.082 | 0.007 |
| | | −0.015 | 0.084 | 0.007 | −0.043 | 0.090 | 0.010 | −0.064 | 0.099 | 0.014 | −0.079 | 0.105 | 0.017 |
| | 1024 | −0.006 | 0.024 | 0.001 | −0.003 | 0.025 | 0.001 | 0.007 | 0.031 | 0.001 | 0.026 | 0.031 | 0.002 |
| | | −0.015 | 0.030 | 0.001 | −0.014 | 0.035 | 0.001 | 0.002 | 0.040 | 0.002 | 0.032 | 0.045 | 0.003 |
| | | −0.002 | 0.042 | 0.002 | −0.003 | 0.042 | 0.002 | −0.005 | 0.042 | 0.002 | −0.005 | 0.042 | 0.002 |
| | | −0.022 | 0.045 | 0.003 | −0.039 | 0.054 | 0.004 | −0.046 | 0.059 | 0.006 | −0.051 | 0.066 | 0.007 |
| $\frac{\pi}{2}$ | 256 | −0.020 | 0.055 | 0.003 | −0.035 | 0.059 | 0.005 | −0.041 | 0.064 | 0.006 | −0.040 | 0.070 | 0.006 |
| | | −0.010 | 0.046 | 0.002 | −0.020 | 0.053 | 0.003 | −0.004 | 0.062 | 0.004 | 0.043 | 0.059 | 0.005 |
| | | 0.002 | 0.094 | 0.009 | 0.000 | 0.094 | 0.009 | 0.000 | 0.093 | 0.009 | 0.005 | 0.084 | 0.007 |
| | | −0.050 | 0.098 | 0.012 | −0.083 | 0.121 | 0.022 | −0.100 | 0.156 | 0.034 | −0.083 | 0.182 | 0.040 |
| | 1024 | −0.012 | 0.022 | 0.001 | −0.015 | 0.024 | 0.001 | −0.007 | 0.028 | 0.001 | 0.014 | 0.034 | 0.001 |
| | | −0.014 | 0.018 | 0.001 | −0.017 | 0.020 | 0.001 | 0.003 | 0.024 | 0.001 | 0.044 | 0.035 | 0.003 |
| | | −0.002 | 0.038 | 0.001 | −0.004 | 0.038 | 0.001 | −0.006 | 0.038 | 0.001 | −0.007 | 0.038 | 0.001 |
| | | −0.039 | 0.046 | 0.004 | −0.064 | 0.069 | 0.009 | −0.061 | 0.096 | 0.013 | −0.023 | 0.097 | 0.010 |

The first row in each cell corresponds to the estimator $\breve{\alpha}(\lambda^0)$, whereas the second, third and fourth correspond to the estimators $\breve{\alpha}(\hat{\lambda}^0)$, $\hat{\alpha}_{\text{LOG}}(\lambda^0)$ and $\hat{\alpha}_{\text{LOG}}(\tilde{\lambda}^0)$, respectively.

**5. Auxiliary results and proofs.** We begin with the proof of Theorem 3.2.

5.1. *Proof of Theorem* 3.2. Let $t = -[vk^{1/2}]$. We examine the case $t > 0$; that for $t < 0$ is similarly handled. First, since Condition C.5 and Lemma 6.10 imply that $|\bar{h}_\psi - h_\psi| = O(k^{-1})$, we have by the definition of $\hat{\alpha}(\lambda_q)$ that

$$\hat{\xi}_n(t/k^{1/2}) = 2^{-1}h_\psi^{-1}\left(\sum_{i=1}^{6}\hat{\xi}_n^{(i)}(t)\right)(1+O_p(k^{-1}))$$

after observing that $\hat{\xi}_n(t/k^{1/2}) = \hat{\xi}_n(v)$, and where

$$\hat{\xi}_n^{(1)}(t) = -\alpha\sum_{p=1}^{k}\psi_p \log(|p-t|_+/p),$$

$$\hat{\xi}_n^{(2)}(t) = -\alpha\sum_{p=1}^{k}\psi_p \log((p+t)/p),$$

$$\hat{\xi}_n^{(3)}(t) = \sum_{p=1}^{k}\psi_p \log\left(\frac{\tilde{f}_{p+s-t}\lambda_{|p-t|_+}^\alpha}{\tilde{f}_{p+s}\lambda_p^\alpha}\right),$$

$$\hat{\xi}_n^{(4)}(t) = \sum_{p=1}^{k}\psi_p \log\left(\frac{\tilde{f}_{s-p-t}\lambda_{p+t}^\alpha}{\tilde{f}_{s-p}\lambda_p^\alpha}\right),$$

$$\hat{\xi}_n^{(5)}(t) = \sum_{p=1}^{k}\psi_p \log\left(\frac{\tilde{f}_{p+s-t}^{-1}\hat{f}_{p+s-t}}{\tilde{f}_{p+s}^{-1}\hat{f}_{p+s}}\right),$$

$$\hat{\xi}_n^{(6)}(t) = \sum_{p=1}^{k}\psi_p \log\left(\frac{\tilde{f}_{s-p-t}^{-1}\hat{f}_{s-p-t}}{\tilde{f}_{s-p}^{-1}\hat{f}_{s-p}}\right),$$

where $\tilde{f}_\ell = (2k_1+1)^{-1}\sum_{j=-k_1}^{k_1} f_{(j+\ell)\mathcal{I}(j+\ell\neq s)+(s+1)\mathcal{I}(j+\ell=s)}$ and $|q|_+ = \max(|q|,1)$.

We examine the behavior of $\hat{\xi}_n^{(i)}(t)$, for $i=1,\ldots,6$, in four propositions. Specifically, Propositions 5.1 and 5.2 deal with the limiting bias of $\sum_{i=1}^{4}\hat{\xi}_n^{(i)}(t)$, although for the proof of Proposition 5.1 we will allow $t < \rho k$ for $0 < \rho < 1$. Proposition 5.3 examines the finite-dimensional limiting distribution of $\hat{\xi}_n^{(5)}(t) + \hat{\xi}_n^{(6)}(t)$ and Proposition 5.4 its tightness. Propositions 5.1–5.4 imply Theorem 3.2.

PROPOSITION 5.1. $\hat{\xi}_n^{(1)}(t) + \hat{\xi}_n^{(2)}(t) = -\bar{\psi}''\alpha\frac{t^2}{k} + O(\frac{t}{k} + \frac{t^{5/2}}{k^{3/2}}).$

PROOF. We only examine $\hat{\xi}_n^{(1)}(t)$, $\hat{\xi}_n^{(2)}(t)$ being identically handled. Assume $\rho < 1/2$ first, so that $0 < t < k/2$. Then

$$(5.1) \quad \hat{\xi}_n^{(1)}(t) = -\alpha \sum_{p=1}^{t} \psi_p \log(|p-t|_+/p) - \alpha \sum_{p=t+1}^{k} \psi_p \log((p-t)/p),$$

where the first term on the right-hand side is $O(k^{-2}t^3)$ because $|\psi_p| \leq Dp^2/k^2$ by Condition C.5 and the integrability of $|u^2 \log((1-u)/u)|$. Next, the second term on the right-hand side of (5.1) is

$$-\alpha \sum_{p=t+1}^{2t} \psi_p \log((p-t)/k)$$

$$-\alpha \sum_{p=t+1}^{k-t} (\psi_{p+t} - \psi_p) \log(p/k) + \alpha \sum_{p=k-t+1}^{k} \psi_p \log(p/k).$$

Proceeding as with the first term on the right-hand side of (5.1), the first term of the last displayed expression is $O(k^{-2}t^3 \log(k/t))$, whereas the last term is $O(k^{-2}t^3)$ by Taylor expansion of $\psi(u)$ around $u = 1$, noting that Condition C.5 implies that $\psi(1) = 0$ and that $\sum_{p=k-t+1}^{k} |\log(p/k)| = O(t^2/k)$ by Taylor expansion of $\log(x)$ around $x = 1$. Finally, the second term of the last displayed expression is

$$(5.2) \quad -\alpha \frac{t}{k} \sum_{p=t+1}^{k-t} \psi'_p \log(p/k) - \frac{\alpha}{2} \frac{t^2}{k} \frac{1}{k} \sum_{p=t+1}^{k-t} \psi''_p \log(p/k) + O\left(\frac{t^{5/2}}{k^{3/2}}\right)$$

by integrability of $|\log u|$ and the fact that $\psi''(u)$ is Lipschitz continuous of order $1/2$ by Condition C.5. By Lemma 6.10 and Condition C.5, the first term of (5.2) is

$$-\alpha t \int_0^1 \psi'(u)(\log u) \, du + O\left(\frac{t}{k}\right) + O\left(\frac{t}{k}\left(\left\{\sum_{p=1}^{t} + \sum_{p=k-t+1}^{k}\right\} |\psi'_p \log(p/k)|\right)\right)$$

$$= \alpha t \int_0^1 \psi(u) u^{-1} \, du + O(k^{-1}t(1 + t^2 k^{-1} \log(k/t))),$$

noting that Condition C.5 implies that $\psi(u)(\log u)|_0^1 = 0$ and $|\psi'_p| \leq Dp/k$ and then proceeding as above. On the other hand, the second term of (5.2) is

$$-\frac{\alpha}{2} \frac{t^2}{k} \int_0^1 \psi''(u)(\log u) \, du + O\left(\frac{t^3}{k^2} \log\left(\frac{k}{t}\right)\right) = -\alpha \frac{\bar{\psi}'' t^2}{2k} + O\left(\frac{t^3}{k^2} \log\left(\frac{k}{t}\right)\right),$$

so that, because $|k^{-1/2} t^{1/2} \log(k/t)| \leq D$ for $0 < t \leq k$, we conclude that

$$\hat{\xi}_n^{(1)}(t) = \alpha t \int_0^1 \psi(u) u^{-1} \, du - \alpha \frac{\bar{\psi}'' t^2}{2k} + O\left(\frac{t}{k} + \frac{t^{5/2}}{k^{3/2}}\right).$$

Now, when $1/2 \leq \rho < 1$, so that $k/2 \leq t < k$, the proof is identical since in this case the left-hand side of (5.1) is $-\alpha \sum_{p=1}^{t} \psi_p \log(|t-p|_+/k) - \alpha \sum_{p=1}^{k-t}(\psi_{p+t} - \psi_p)\log(p/k) + \alpha \sum_{p=k-t+1}^{k} \psi_p \log(p/k)$. Then proceed as above. Proceeding similarly, $\hat{\xi}_n^{(2)}(t) = -\alpha t \int_0^1 \psi(u) u^{-1} \, du - \alpha \frac{\bar{\psi}'' t^2}{2k} + O(\frac{t}{k} + \frac{t^{5/2}}{k^{3/2}})$. From here the conclusion is obvious. $\square$

PROPOSITION 5.2. $\hat{\xi}_n^{(3)}(t) + \hat{\xi}_n^{(4)}(t) = o(1)$.

PROOF. We only examine $\hat{\xi}_n^{(3)}(t)$, as $\hat{\xi}_n^{(4)}(t)$ is similar. By definition $\hat{\xi}_n^{(3)}(t)$ is

$$-\sum_{p=1}^{2k_1} \psi_p a_p - \sum_{p=2k_1+1}^{k} \psi_p a_p - \sum_{p=1}^{k} \psi_p \tilde{g}_p, \tag{5.3}$$

where

$$a_p = \log(f_{(p+s)\mathcal{I}(p\neq 0)+(s+1)\mathcal{I}(p=0)}^{-1} \tilde{f}_{p+s})$$
$$- \log(f_{(p+s-t)\mathcal{I}(p\neq t)+(s+1)\mathcal{I}(p=t)}^{-1} \tilde{f}_{p+s-t})$$

and

$$\tilde{g}_p = \log(\lambda_p^{\alpha} f_{(p+s)\mathcal{I}(p\neq 0)+(s+1)\mathcal{I}(p=0)})$$
$$- \log(\lambda_{|p-t|_+}^{\alpha} f_{(p+s-t)\mathcal{I}(p\neq t)+(s+1)\mathcal{I}(p=t)}).$$

Since by Condition C.1 and $|\lambda^0 - \lambda_s| \leq \frac{\pi}{n}$, $D^{-1}(k_1/p)^{\alpha} < \lambda_{k_1}^{\alpha} f_{p+s} < D(k_1/p)^{\alpha}$, it implies that, for $p \leq 2k_1$, $|a_p| = O(\log(k_1/p))$ by Lemma 6.1. Note that $a_t = O(\log k_1)$. Hence the absolute value of the first term of (5.3) is bounded by

$$D|\psi_t|\log(k_1) + D \sum_{p=1; p\neq t}^{2k_1} |\psi_p| \log\left(\frac{k_1}{p}\right) = O\left(\frac{k_1^3}{k^2}\right) = o(1)$$

by Condition C.4 and because Condition C.5 implies that $|\psi_p| \leq D(p/k)^2$. The absolute value of the second term of (5.3) is bounded by

$$D \sum_{p=2k_1+1}^{k_1[\log^{1/3} k_1]} |\psi_p a_p| + D \sum_{p=k_1[\log^{1/3} k_1]+1}^{k} |\psi_p a_p|$$
$$= O\left(\frac{k_1^3 \log k_1}{k^2} + \frac{t k_1^2}{k^2}\left(1 + \log\left(\frac{k}{t}\right)\right)\right) = o(1),$$

where for the first term on the left-hand side we have used the fact that by Lemma 6.1(a), $D^{-1} < |f_{p+s}^{-1} \tilde{f}_{p+s}| < D$ and then Condition C.5, and for the

second term on the left-hand side the fact that by Lemma 6.1(a), $|f_{p+s}^{-1}\tilde{f}_{p+s} - 1| = O(p^{-2}k_1^2)$, which implies that $|\log(f_{p+s}^{-1}\tilde{f}_{p+s})| = O(p^{-2}k_1^2)$ by the mean value theorem, and then Lemma 6.4 with $\nu_p = k_1^{-2}\log(f_{p+s}^{-1}\tilde{f}_{p+s})$ there.

To complete the proof, it remains to show that the third term of (5.3) is $o(1)$. By Condition C.1, $\sum_{p=1}^{k}\psi_p\tilde{g}_p$ is

$$
(5.4) \quad \begin{aligned} &\sum_{p=1}^{k}\psi_p(\log(g(\lambda_{p+s})) - \log(g(\lambda_{p+s-t}))) \\ &- \alpha\sum_{p=1}^{k}\psi_p(\log(|\lambda_{p+s} - \lambda^0|\lambda_p^{-1}) - \log(|\lambda_{p-t+s} - \lambda^0|_+\lambda_{|p-t|_+}^{-1})). \end{aligned}
$$

Denote the first and second derivatives of $\log(g(\lambda))$ by $h(\lambda)$ and $h'(\lambda)$, respectively. The first term of (5.4) is

$$
-\left(\frac{2\pi t}{n}\right)\sum_{p=1}^{k}\psi_p h(\lambda_{p+s}) - \left(\frac{2\pi t}{n}\right)^2 \frac{1}{2}\sum_{p=1}^{k}\psi_p h'(\lambda_{p+s-\theta(p)t})
$$

$$
= -\frac{4\pi^2 t}{n^2}h'(\lambda_s)\sum_{p=1}^{k}p\psi_p + o(1),
$$

where $\theta(p) \in (0,1)$, by Condition C.4 and because Lemma 6.10 and Condition C.5 imply that $k^{-1}\sum_{p=1}^{k}\psi_p = O(k^{-1})$, so that

$$
\sum_{p=1}^{k}\psi_p h(\lambda_{p+s}) = \sum_{p=1}^{k}\psi_p(h(\lambda_{p+s}) - h(\lambda_s)) + O(1)
$$

$$
= \frac{2\pi}{n}h'(\lambda_s)\sum_{p=1}^{k}p\psi_p + O(k^3 n^{-2} + 1).
$$

Next, the second term of (5.4) is

$$
-\alpha\left\{\sum_{p=1}^{2t} + \sum_{p=2t+1}^{k}\right\}\psi_p(\log((\lambda_{p+s} - \lambda^0)\lambda_p^{-1}) - \log(|\lambda_{p-t+s} - \lambda^0|_+\lambda_{|p-t|_+}^{-1})).
$$

The contribution due to $\sum_{p=1}^{2t}$ is $o(1)$ by Condition C.5 and then Condition C.4, whereas the contribution due to $\sum_{p=2t+1}^{k}$, by Taylor expansion of $\log(x)$, is

$$
-\alpha n\left(\frac{\lambda_s - \lambda^0}{2\pi}\right)\sum_{p=2t+1}^{k}\psi_p\left(\frac{1}{p} - \frac{1}{p-t}\right)
$$

$$+\frac{\alpha}{2}n^2\left(\frac{\lambda_s-\lambda^0}{2\pi}\right)^2\sum_{p=2t+1}^{k}\psi_p\left\{\frac{1}{p^2(s)}-\frac{1}{p^2(s-t)}\right\}$$

$$= O(tk^{-1}\log k) + O(k^{-1} + tk^{-2}\log k) = o(1),$$

where $p(s-a)$ is an intermediate point between $p-a$ and $(p-a)+n(\lambda_s-\lambda^0)/(2\pi)$, and then because $n|\lambda^0 - \lambda_s| \leq \pi$, $|p^{-j}\psi_p| \leq Dk^{-j}$ for $j=1,2$, by Condition C.5 and the fact that $|p^2 p^{-2}(s)| + |(p-t)^2 p^{-2}(s-t)| \leq D$. So, we conclude that

$$\hat{\xi}_n^{(3)}(t) = -\frac{4\pi^2 t}{n^2}h'(\lambda_s)\sum_{p=1}^{k}p\psi_p + o(1).$$

Similarly, $\hat{\xi}_n^{(4)}(t) = 4\pi^2 t h'(\lambda_s)\sum_{p=1}^{k} p\psi_p/n^2 + o(1)$. From here the conclusion of the proposition is obvious. □

PROPOSITION 5.3. *The finite-dimensional distributions of $\hat{\xi}_n^{(5)}(t)+\hat{\xi}_n^{(6)}(t)$ converge to those of a normal random variable.*

PROOF. By the Wold device, it suffices to show that, for any finite $l > 0$,

$$\sum_{i=1}^{l}\phi_i(\hat{\xi}_n^{(5)}(t_i) + \hat{\xi}_n^{(6)}(t_i)) \xrightarrow{d} N\left(0, \varsigma \sum_{i,j=1}^{l}\phi_i\phi_j v_i v_j\right),$$

where $\phi_i$ satisfies $\sum_{i=1}^{l}|\phi_i|^2 = 1$ and $v_i = \lim t_i/k^{1/2}$. Denoting $\tilde{f}_\ell^{-1}\hat{f}_\ell$ by $\hat{g}_\ell$,

$$\hat{\xi}_n^{(5)}(t) = -\sum_{p=1}^{2t}\psi_p\log(\hat{g}_{p+s}/\hat{g}_{p+s-t}) - \sum_{p=2t+1}^{k}\psi_p\log(\hat{g}_{p+s}/\hat{g}_{p+s-t})$$

(5.5)
$$:= b_{1t} + b_{2t}.$$

We begin by showing that $b_{1t} = o_p(1)$. In particular, we will show something stronger than needed, that is, that for each $\varepsilon > 0$ there exists $n_0$ such that

(5.6) $$\Pr\left\{\sup_{t_1 \leq q \leq t_2}|b_{1q} - b_{1t_1}| > \varepsilon\right\} < DM^3 k^{-1/2}\varepsilon^{-1}\log^{-1} k_1$$

for all $n \geq n_0$ and $0 \leq t_1 < t_2 \leq [k^{1/2}M]$. Because $t_2 \leq [k^{1/2}M]$,

$$\sup_{t_1 \leq q \leq t_2}|b_{1q} - b_{1t_1}| \leq 2\sup_{0\leq|q|\leq 2[k^{1/2}M]}|\log\hat{g}_{q+s}|\sum_{p=1}^{2[k^{1/2}M]}|\psi_p|$$

$$\leq D\sup_{0\leq|q|\leq 2[k^{1/2}M]}|\log\hat{g}_{q+s}|M^3 k^{-1/2},$$

since Condition C.5 implies that $\sum_{p=1}^{a} |\psi_p| \leq Da^3/k^2$. Then by Markov's inequality the left-hand side of (5.6) is bounded by

$$DM^3 k^{-1/2} \varepsilon^{-1} E \sup_{0 \leq |q| \leq 2[k^{1/2}M]} |\log \hat{g}_{q+s}| \leq DM^3 k^{-1/2} \varepsilon^{-1} \log^{-1} k_1,$$

because Lemma 6.1 and the definition of $\hat{f}_{q+s}$ imply that $\hat{g}_{q+s}$ is bounded from below by $D^{-1} n^{-h}$, for some $h \geq 1$, and for $x > D^{-1} n^{-h}$,

(5.7) $$|\log x - (x-1)| \leq D(x-1)^2 \log n,$$

by Condition C.4, $|\log n / \log k_1| \leq D$, and the fact that by Lemma 6.3 and Condition C.4, for some $\beta > 0$,

(5.8) $$E\left(\sup_{q:\,|q| \leq 2k_1} |\hat{g}_{q+s} - 1|^\mu\right) \leq \frac{D}{\log^{1+\mu} k_1},$$

(5.9) $$E\left(\sup_{q:\,2k_1 < |q|} |\hat{g}_{q+s} - 1|^\mu\right) \leq \frac{D}{k_1^{\beta\mu}}.$$

So (5.6) holds, which implies that $b_{1t} = o_p(1)$.

Next we look at the second term on the right-hand side of (5.5), that is, $b_{2t}$. Denoting

(5.10) $$b_{3t} = \sum_{p=2t+1}^{k-t} (\psi_{p+t} - \psi_p) \log(\hat{g}_{p+s}),$$

we have that $b_{2t} - b_{3t}$ is

(5.11) $$\sum_{p=t+1}^{2t} \psi_{p+t} \log(\hat{g}_{p+s}) - \sum_{p=k-t+1}^{k} \psi_p \log(\hat{g}_{p+s}).$$

Because by Condition C.5, $\{\sum_{p=t+1}^{2t} |\psi_{p+t}| + \sum_{p=q+1}^{2q} |\psi_{p+q}|\} \leq Dk^{-2}(t^3 + q^3) \leq DM^3 k^{-1/2}$, the first term of (5.11) satisfies

(5.12) 
$$\Pr\left\{\sup_{t_1 \leq q \leq t_2} \left|\sum_{p=t_1+1}^{2t_1} \psi_{p+t_1} \log(\hat{g}_{p+s}) - \sum_{p=q+1}^{2q} \psi_{p+q} \log(\hat{g}_{p+s})\right| > \varepsilon\right\}$$
$$\leq \frac{DM^3}{\varepsilon k^{1/2} \log k_1}$$

by Markov's inequality and (5.7)–(5.8), whereas the second term satisfies

$$\Pr\left\{\sup_{t_1 \leq q \leq t_2} \left|\left(\sum_{p=k-t_1+1}^{k} - \sum_{p=k-q+1}^{k}\right) \psi_p \log(\hat{g}_{p+s})\right| > \varepsilon\right\} \leq \frac{DM}{\varepsilon \log k_1} \left(\frac{t_2 - t_1}{k^{1/2}}\right)$$

by (5.7) and (5.9) and because $\sum_{p=k-q+1}^{k} - \sum_{p=k-t_1+1}^{k} = \sum_{p=k-q+1}^{k-t_1}$ and by Condition C.5, $\sup_{t_1 \leq q \leq t_2} \sum_{p=k-q+1}^{k-t_1} |\psi_p| \leq \sum_{p=k-t_2}^{k-t_1} |\psi_p| \leq Dk^{-1}(t_2^2 - t_1^2) \leq DMk^{-1/2}(t_2 - t_1)$.

So, (5.6), (5.12) and the last inequality imply that for any $\varepsilon > 0$ there exists $n_0$ such that for all $n \geq n_0$,

$$\Pr\left\{\sup_{t_1 \leq q \leq t_2} |(\hat{\xi}_n^{(5)}(q) - b_{3q}) - (\hat{\xi}_n^{(5)}(t_1) - b_{3t_1})| > \varepsilon\right\}$$
(5.13)
$$\leq \frac{DM}{\varepsilon k^{1/2} \log k_1}(M^2 + (t_2 - t_1)).$$

Clearly (5.13) implies that $\sup_{0 \leq t \leq 2[k^{1/2}M]} |\hat{\xi}_n^{(5)}(t) - b_{3t}| = o_p(1)$. Proceeding similarly, we have that for any $\varepsilon > 0$ there exists $n_0$ such that for all $n \geq n_0$,

$$\Pr\left\{\sup_{t_1 \leq q \leq t_2} |(\hat{\xi}_n^{(6)}(q) - b_{4q}) - (\hat{\xi}_n^{(6)}(t_1) - b_{4t_1})| > \varepsilon\right\}$$
(5.14)
$$\leq \frac{DM}{\varepsilon k^{1/2} \log k_1}(M^2 + (t_2 - t_1)),$$

where

(5.15) $$b_{4t} = \sum_{p=3t+1}^{k} (\psi_{p-t} - \psi_p) \log(\hat{g}_{s-p}).$$

Next we examine $b_{3t}$ and $b_{4t}$. Denoting $\vartheta_\ell = \hat{g}_\ell - 1$ and writing

(5.16) $$\tilde{b}_t = \sum_{p=2t+1}^{k-t} (\psi_{p+t} - \psi_p)\vartheta_{p+s}, \qquad \tilde{\tilde{b}}_t = \sum_{p=3t+1}^{k} (\psi_{p-t} - \psi_p)\vartheta_{s-p},$$

Lemma 6.5 implies that $b_{3t} = \tilde{b}_t + t/k^{1/2}o_p(1)$ and $b_{4t} = \tilde{\tilde{b}}_t + t/k^{1/2}o_p(1)$, where the $o_p(1)$ term is uniformly in $t \leq \rho k$, for $\rho < 1/3$.

So it remains to examine $\tilde{b}_t$ and $\tilde{\tilde{b}}_t$. By Taylor expansion of $\psi_p$,

(5.17) $$\tilde{b}_t = \frac{t}{k} \sum_{p=2t+1}^{k-t} \psi'_p \vartheta_{p+s} + \frac{1}{2}\frac{t^2}{k^2} \sum_{p=2t+1}^{k-t} \psi''\left(\frac{p}{k} + \delta\frac{t}{k}\right)\vartheta_{p+s},$$

where $\delta = \delta(t) \in (0, 1)$. The first term on the right-hand side of (5.17) is

$$\frac{t}{k}\sum_{p=2k_1+1}^{k-2k_1} \psi'_p \vartheta_{p+s} + \frac{t}{k}\left\{\sum_{p=2t+1}^{2k_1} + \sum_{p=k-2k_1+1}^{k-t}\right\}\psi'_p \vartheta_{p+s}.$$

By Lemma 6.2 and Conditions C.4 and C.5, the second term of the last displayed expression is clearly $tk_1^{1/2}/kO_p(1 + k^{-1}k_1^{\alpha+1/2}\mathcal{I}(\alpha \geq 1/2)) = t/k^{1/2}o_p(1)$, where the $o_p(1)$ term does not depend on $t \leq \rho k$. On the other hand, writing $\eta_j = f_{j+s}^{-1} I_{j+s} - 1$, the first term is

$$
\begin{aligned}
& \frac{t}{k} \sum_{p=2k_1+1}^{k-2k_1} \frac{\psi'_p}{2k_1+1} \sum_{j=-k_1}^{k_1} \eta_{p+j} \\
& \quad + \frac{t}{k} \sum_{p=2k_1+1}^{k-2k_1} \frac{\psi'_p}{2k_1+1} \sum_{j=-k_1}^{k_1} \left( \frac{f_{p+j+s}}{\tilde{f}_{p+s}} - 1 \right) \eta_{p+j}.
\end{aligned}
\tag{5.18}
$$

Since $|\tilde{f}_{p+s}^{-1} f_{p+s} - 1| = \tilde{f}_{p+s}^{-1} f_{p+s} |1 - f_{p+s}^{-1} \tilde{f}_{p+s}| = O(k_1^2/p^2)$ by Lemma 6.1, $|f_{p+s}^{-1} f_{p+j+s} - 1| \leq D \frac{k_1}{p} \frac{|j|}{k_1}$ by Condition C.3, and by an obvious extension of Robinson [32], $E|\sum_{j=-k_1}^{k_1} c_j \eta_{p+j}| = O(k_1^{1/2})$ for any $|c_j = c(j/k_1)| \leq D$, we obtain that the first absolute moment of the second term of (5.18) is by Condition C.5 and then Condition C.4, $tk^{-1/2} O(k_1^{1/2} k^{-1/2} \sum_{p=1}^{k} |\psi'_p|/p) = t/k^{1/2} o(1)$, where the $o(1)$ term does not depend on $t \leq \rho k$.

On the other hand, after rearranging subindices, the first term of (5.18) is

$$
\begin{aligned}
& \frac{t}{k} \sum_{p=2k_1+1}^{k-2k_1} \eta_p \left( \frac{1}{2k_1+1} \sum_{j=1}^{2k_1} \psi'_{p+j} \right) \\
& \quad + \frac{t}{k} \sum_{p=k_1+1}^{2k_1} \eta_p \left( \frac{1}{2k_1+1} \sum_{j=1}^{p-k_1} \psi'_{j+k_1} \right) \\
& \quad + \frac{t}{k} \sum_{p=k-2k_1+1}^{k-k_1} \eta_p \left( \frac{1}{2k_1+1} \sum_{j=p}^{k-k_1} \psi'_{j-k_1} \right).
\end{aligned}
$$

After standard calculations and routine application of Robinson's [32] Theorem 2, the last two terms are $t/k O_p(k_1^{1/2}) = t/k^{1/2} o_p(1)$ by Conditions C.4 and C.5, whereas the first term of the last displayed expression is

$$
\frac{t}{k} \sum_{p=2k_1+1}^{k-2k_1} \psi'_p \eta_p + \frac{t}{k} \sum_{p=2k_1+1}^{k-2k_1} \psi'_p \eta_p \left( \frac{1}{\psi'_p} \left[ \frac{1}{2k_1+1} \sum_{j=1}^{2k_1} \psi'_{p+j} \right] - 1 \right).
$$

We note that $\psi'(u)$ continuous by Condition C.5 implies that $\psi'_{p+j}/\psi'_p \to 1$ as $k_1/p \to 0$, and hence the expression inside the parentheses converges to zero as $p \to \infty$. So by Toeplitz's lemma we conclude that the last displayed expression, and therefore also the first term of (5.17), is

$$
\frac{t}{k^{1/2}} \left( \frac{1}{k^{1/2}} \sum_{p=2k_1+1}^{k-2k_1} \psi'_p \eta_p + o_p(1) \right).
$$

Proceeding similarly as with the first term on the right-hand side of (5.17), the second term of (5.17) is $k^{-3/2}t^2(k^{-1/2}\sum_{p=2k_1+1}^{k-2k_1}\psi''(\frac{p}{k}+\delta\frac{t}{k})\eta_p + o_p(1))$, so that

$$(5.19) \qquad \tilde{b}_t = \frac{t}{k^{1/2}}\left(\frac{1}{k^{1/2}}\sum_{p=2k_1+1}^{k-2k_1}\psi'_p\eta_p + \frac{t}{k}O_p(1) + o_p(1)\right),$$

where the $o_p(1)$ and $O_p(1)$ terms are uniform in $t \le \rho k$. Similarly, we obtain that

$$(5.20) \qquad \tilde{\tilde{b}}_t = \frac{t}{k^{1/2}}\left(\frac{1}{k^{1/2}}\sum_{p=2k_1+1}^{k-2k_1}\psi'_p\eta_{-p} + \frac{t}{k}O_p(1) + o_p(1)\right),$$

where the $o_p(1)$ and $O_p(1)$ terms are uniform in $t \le \rho k$.

Thus, (5.13), (5.14), Lemma 6.5, (5.19) and (5.20) imply that, for $t_i \le [Mk^{1/2}]$,

$$(5.21) \quad \sum_{i=1}^{l}\phi_i(\hat{\xi}_n^{(5)}(t_i)+\hat{\xi}_n^{(6)}(t_i)) = \sum_{i=1}^{l}\phi_i(\tilde{b}_{t_i}+\tilde{\tilde{b}}_{t_i}) + \frac{t_i}{k^{1/2}}o_p(1)$$
$$\xrightarrow{d} N\left(0, \varsigma\sum_{i,j=1}^{l}\phi_i\phi_j v_i v_j\right)$$

by Robinson's [32] Theorem 2 and Toeplitz's lemma, since $|k^{-1}\sum_{\ell=t+1}^{k-t}(\psi'_\ell)^2 - \varsigma| = o(1)$ by Lemma 6.10. $\square$

PROPOSITION 5.4. $\hat{\xi}_n^{(5)}(t)$ and $\hat{\xi}_n^{(6)}(t)$ are tight.

PROOF. Write $c_t = \hat{\xi}_n^{(5)}(t) - b_{3t}$, where $b_{3t}$ is given in (5.10). To show that $\hat{\xi}_n^{(5)}(t)$ is tight it suffices to show that $c_t$ and $b_{3t}$ are tight. Since the finite-dimensional distributions of $c_t$ converge to zero [cf. (5.13)], Billingsley's [4] Theorem 15.4 implies that $c_t$ is tight if for each $\varepsilon > 0$ and $\nu > 0$ there exists a $\delta \in (0,1)$ such that

$$(5.22) \qquad \Pr\{\vartheta''(c_t,\delta) \ge \varepsilon\} \le \nu$$

holds for all $n \ge n_0$, where

$$\vartheta''(c_t,\delta) = \sup\min\{|c_t - c_{t_1}|, |c_{t_2} - c_t|\},$$

and the supremum is over $t_1, t$ and $t_2$ satisfying $t_1 \le t \le t_2$ with $t_2 - t_1 \le \delta[k^{1/2}M]$ and $\delta \in (0,1)$. Observe that we can assume $k^{-1/2} \le [t_2/k^{1/2}M] - [t_1/k^{1/2}M]$. If $[t_2/k^{1/2}M] - [t_1/k^{1/2}M] < k^{-1/2}$, then either $t_1$ and $t$ lie in the same subinterval $[(p-1)/M, p/M)$ or else $t$ and $t_2$ do; in either of these cases the left-hand side of (5.22) vanishes.

Inequalities (14.9) and (14.46) in [4] imply that (5.22) holds if
$$\Pr\{\vartheta(c_t,\delta) \geq \varepsilon\} \leq \nu,$$
for some $0 < \delta \leq 1$, where
$$\vartheta(c_t,\delta) = \sup_{|(t-v)/(k^{1/2}M)| < \delta} |c_t - c_v|.$$
(Observe that as $c_t$ converges in probability to zero, which has continuous paths, the Skorohod metric can be replaced by the uniform topology.) By the corollary of Billingsley's [4] Theorem 8.3, it suffices to show that

$$(5.23) \qquad \sum_{i=1}^{r} \Pr\left\{\sup_{t_{i-1} \leq v \leq t_i} |c_v - c_{t_{i-1}}| \geq \varepsilon/3\right\} \leq \nu,$$

where $2^{-1}\delta < [k^{1/2}M]^{-1}(t_i - t_{i-1}) < \delta$ and $0 = t_0 < t_1 < \cdots < t_r = [k^{1/2}M]$. But this is the case since by (5.13),
$$\Pr\left\{\sup_{t_{i-1} \leq v \leq t_i} |c_v - c_{t_{i-1}}| \geq \varepsilon/3\right\} \leq \frac{DM^3\delta}{\varepsilon \log k_1}.$$

Now choose $n_0$ such that $DM^3\varepsilon^{-1}\log^{-1} k_1 < \nu$ since $r \leq 2[\delta^{-1}]$ to obtain (5.23). Proceeding similarly, but using (5.14) instead of (5.13), $\hat{\xi}_n^{(6)}(t) - b_{4t}$ is also tight.

Next we show the tightness condition for $b_{3t}$; the proof for $b_{4t}$ is similar and is omitted. Consider $t < q$. Then $b_{3t} - b_{3q}$ is

$$\sum_{p=2t+1}^{k-t} (\psi_{p+t} - \psi_p) \log \hat{g}_{p+s} - \sum_{p=2q+1}^{k-q} (\psi_{p+q} - \psi_p) \log \hat{g}_{p+s}$$

$$(5.24) \qquad = \left\{\sum_{p=2t+1}^{2q} + \sum_{p=k-q+1}^{k-t}\right\}(\psi_{p+t} - \psi_p) \log \hat{g}_{p+s}$$

$$+ \sum_{p=2q+1}^{k-q} (\psi_{p+t} - \psi_{p+q}) \log \hat{g}_{p+s}.$$

The first term on the right-hand side of (5.24) is tight, as we now show. Because by Condition C.5,

$$|\psi_{p+t} - \psi_p| \leq Dk^{-1}t, \qquad \left|\sum_{p=2t+1}^{2q} 1 + \sum_{p=k-q+1}^{k-t} 1\right| \leq 3(q-t),$$

abbreviating the first term on the right-hand side of (5.24) by $\zeta_{t,q}$, we obtain that
$$\sum_{i=1}^{r} \Pr\left\{\sup_{t_{i-1} \leq q \leq t_i} |\zeta_{t_{i-1},q}| \geq \varepsilon\right\} \leq \sum_{i=1}^{r} \Pr\left\{D\frac{t_i - t_{i-1}}{k^{1/2}} \sup_{p=2t+1,\ldots,k} |\log \hat{g}_{p+s}| \geq \varepsilon\right\}$$

$$\leq \frac{DM}{\varepsilon \log k_1}$$

by Markov's inequality and (5.7)–(5.8) with $\mu = 1$ there. Then choose $n_0$ such that $DM\varepsilon^{-1}\log^{-1} k_1 < \nu$ to complete.

Next, Taylor expansion implies that the second term on the right-hand side of (5.24) is

$$(5.25) \quad \left(\frac{t-q}{k}\right) \sum_{p=2q+1}^{k-q} \psi'_{p+t} \log \hat{g}_{p+s} + \frac{1}{2}\left(\frac{t-q}{k}\right)^2 \sum_{p=2q+1}^{k-q} \psi''_{p+\ell} \log \hat{g}_{p+s},$$

where $\ell$ is an intermediate point between $t$ and $q$.

The second term of (5.25) is tight as we now show. Proceeding as with the proof of tightness of $c_t$, it suffices to show that for all $\nu$ and $\varepsilon > 0$ there exists $n_0$ such that

$$\sum_{i=1}^{r} \Pr\left\{ \sup_{t_{i-1} \leq q \leq t_i} \left| \left(\frac{q-t_{i-1}}{k}\right)^2 \sum_{p=2q+1}^{k-q} |\psi''_{p+\ell} \log \hat{g}_{p+s}| \right| \geq \varepsilon \right\} \leq \nu$$

for all $n \geq n_0$ and $0 \leq t_0 < \cdots < t_r \leq [k^{1/2}M]$. But by (5.7) and (5.8) and the fact that $|\psi''(u)| \leq D$ by Condition C.5, the left-hand side of the last displayed inequality is bounded by

$$D\varepsilon^{-1}k^{-1}\log^{-1} k_1 \sum_{i=1}^{r}(t_i - t_{i-1})^2 \leq DM^2\varepsilon^{-1}\log^{-1} k_1 \sum_{i=1}^{r} \delta^2 < \nu,$$

since $r \leq 2[\delta^{-1}]$ and $2^{-1}\delta < [k^{1/2}M]^{-1}(t_i - t_{i-1}) < \delta$.

To finish the proof it remains to examine the first term of (5.25), denoted by $d_{t,q}$. Since from the proof of Proposition 5.3, the finite-dimensional distributions of $d_{t,q}$ converge to those of the limiting Gaussian process which has continuous paths, by Billingsley's [4] Theorem 15.4, it implies that it suffices to check that

$$\Pr\{\vartheta''(d_{t,q}, \delta) \geq \varepsilon\} \leq \nu$$

for some $n \geq n_0$. Now by Billingsley's [4] Theorem 15.6, it suffices to check the moment condition

$$(5.26) \quad E|d_{t,q}d_{q,v}|^{\beta_1} \leq D \left|\frac{t-v}{k^{1/2}}\right|^{\beta_2}$$

for $t \leq q \leq v$ and some $\beta_1 > 0$ and $\beta_2 > 1$. Write

$$(5.27) \quad d_{t,q} = \left(\frac{q-t}{k}\right) \sum_{p=2q+1}^{k-q} \psi'_{p+t}(\hat{g}_{p+s} - 1)$$

$$+ \left(\frac{q-t}{k}\right)\left\{ \sum_{p=2q+1}^{3k_1} + \sum_{p=3k_1+1}^{k-q} \right\} \psi'_{p+t}(1 - \hat{g}_{p+s} + \log \hat{g}_{p+s}).$$

Using (5.7),ced the first absolute moment of the second term on the right-hand side of (5.27) is bounded by

$$D\left|\frac{q-t}{k}\right|\log n\left\{\sum_{p=2q+1}^{3k_1} + \sum_{p=3k_1+1}^{k-q}\right\}|\psi'_{p+t}|E(\hat{g}_{p+s}-1)^2$$

$$\leq D|q-t|\log n\left(\frac{k_1^{2\alpha}}{k^2}\mathcal{I}(\alpha\geq 1/2) + \frac{1}{k_1}\right)$$

$$\leq D\left|\frac{t-q}{k^{1/2}}\right|^{1+\xi}\log n\left(\frac{k^{(\xi+1)/2}}{k_1} + \frac{k_1^{2\alpha}k^{(1+\xi)/2}}{k^2}\mathcal{I}(\alpha\geq 1/2)\right)$$

$$\leq D\left|\frac{t-q}{k^{1/2}}\right|^{1+\xi}$$

for some $0 < \xi < 3(1-\alpha)$ by Lemma 6.2 and observing that $\sum_{p=2q+1}^{3k_1}|\psi'_{p+t}| \leq Dk_1^2 k^{-1} = o(k^{1/2})$ by Conditions C.4 and C.5. So, the second term on the right-hand side of (5.27) satisfies (5.26), which follows by the Cauchy–Schwarz inequality, choosing $\beta_1 = 1/2$ and $\beta_2 = 1 + \xi$, and the fact that $(q-t)(v-q) \leq (v-t)^2$. Finally, proceeding as with the proof of $\tilde{b}_t$ given in (5.19),

$$E\left|\left(\frac{q-t}{k}\right)\sum_{p=2q+1}^{k-q}\psi'_{p+t}(\hat{g}_{p+s}-1)\right|^2 \leq D\left(\frac{q-t}{k^{1/2}}\right)^2,$$

which implies that the first term on the right-hand side of (5.27) satisfies (5.26) choosing $\beta_1 = 1$ and $\beta_2 = 2$, and noting that $(q-t)(v-q) \leq (v-t)^2$. So, we conclude that $b_{3t}$ is tight. Proceeding similarly, $b_{4t}$ is also tight, which completes the proof. □

5.2. *Proof of Corollary* 3.3. (a) By Theorem 3.2, $\hat{\xi}_n(v) \overset{\text{weakly}}{\Longrightarrow} \xi(v)$ in $\mathbb{D}[-M, M]$ for any arbitrary $M > 0$. Next, since the limiting Gaussian process $\xi(v)$ has continuous paths, that is, it belongs to $\mathbb{C}[-M, M]$, the Skorohod metric can be replaced by the uniform topology. On the other hand, by Eddy [10], the *argmaximum* is a continuous functional in the set of parabolas with fixed second derivatives in $\mathbb{C}[-M, M]$. So, by van der Vaart and Wellner's [34] Theorem 3.2.2, we obtain that

$$\hat{v}_n = \arg\max_v \hat{\xi}_n(v) \overset{d}{\to} \arg\max_v \xi(v) = v^*,$$

where $v^* = \Psi^{1/2}X$ and $X = N(0,1)$. Observe that Theorem 3.1 shows that $\Pr\{|\hat{v}_n| < L\} > 1 - \delta$ for $n$ sufficiently large. This together with Problem 1.3.9 in [34], page 27, implies that $\hat{v}_n$ is uniformly tight.

But by construction, $\hat{\lambda}^0 = \lambda^0 + n^{-1}(2\pi k^{1/2})\hat{v}_n + O(n^{-1})$, that is, (3.2), and hence

$$(2\pi k^{1/2})^{-1} n(\hat{\lambda}^0 - \lambda^0) = \hat{v}_n + O(k^{-1/2}) \xrightarrow{d} v^* = \Psi^{1/2} X.$$

(b) As in (a), the limit process is $\xi(v)$, where from the definition of $\hat{\lambda}^0$, $v \geq 0$. Thus, if $X$ takes a positive value, the restriction $v \geq 0$ is not binding and the maximum of $\xi(v)$ is achieved at $v^*$. However, when $X$ takes a negative value, the restriction is binding and thus the maximum is at $v = 0$ due to the parabolic structure of $\xi(v)$.

(c) The proof is identical to part (b) once the wording positive (negative) is replaced by negative (positive). □

5.3. *Proof of Theorem* 3.1. Similarly to the proof of Theorem 3.2, it suffices to show the theorem with $\bar{h}_\psi$ replaced by $h_\psi$. With that replacement and recalling that $\hat{g}_p = \tilde{f}_p^{-1} \hat{f}_p$, $\hat{\alpha}(\lambda_q)$ becomes

$$\frac{h_\psi^{-1}}{2k} \sum_{p=1}^{k} \psi_p (\log \hat{g}_{p+q} + \log \hat{g}_{q-p})$$

$$+ \frac{h_\psi^{-1}}{2k} \sum_{p=1}^{k} \psi_p \log(\lambda_{|p+q-s|_+}^\alpha \tilde{f}_{p+q})$$

(5.28)

$$+ \frac{h_\psi^{-1}}{2k} \sum_{p=1}^{k} \psi_p \log(\lambda_{|q-p-s|_+}^\alpha \tilde{f}_{q-p})$$

$$- \alpha \frac{h_\psi^{-1}}{2k} \sum_{p=1}^{k} \psi_p (\log(2\pi|p+q-s|_+/n) + \log(2\pi|q-p-s|_+/n)).$$

Because (5.7)–(5.9) imply that $\sup_{\ell=0,\ldots,[n/2]} |\log \hat{g}_\ell| = o_p(1)$, we then have that, uniformly in $q$, the first term of (5.28) converges to zero in probability. Next, consider the second term of (5.28). (Recall that $0 \leq q \leq [n/2]$.) If $q < \max\{0, s - [k_1 \log k_1] - k\}$ or $s + [k_1 \log k_1] < q$, this is

(5.29)

$$\frac{h_\psi^{-1}}{2k} \sum_{p=1}^{k} \psi_p \left\{ \log\left(\frac{\tilde{f}_{p+q}}{f_{p+q}}\right) - \log\left(\frac{g(\lambda_q)}{g(\lambda_{p+q})}\right) \right.$$
$$\left. + \log\left(\frac{\lambda_{p+q-s}^\alpha}{|\lambda_{p+q} - \lambda^0|^\alpha}\right) \right\} + o(1),$$

because Lemma 6.10 and Condition C.5 imply that $k^{-1} \sum_{p=1}^{k} \psi_p = O(k^{-1})$ and by Condition C.1, $f_p = |\lambda_p - \lambda^0|^{-\alpha} g(\lambda_p)$. But by Lemma 6.1(a) and

Taylor expansion of $\log(z)$ around $z = 1$, (5.29) is bounded in absolute value by

$$\frac{D}{k}\sum_{p=1}^{k}|\psi_p|\left(\left|\log\left(\frac{g(\lambda_q)}{g(\lambda_{p+q})}\right)\right| + \frac{k_1^2}{|q+p-s|^2} + \alpha\log\left(\frac{\lambda|_{p+q-s|}}{|\lambda_{p+q} - \lambda^0|}\right)\right) + o(1).$$

The contribution due to the second term inside the brackets is easily shown to be $o(1)$, as is the contribution due to the third term by Taylor expansion and the fact that $n|\lambda_s - \lambda^0| < \pi$. Finally, by the mean value theorem, the first term inside the brackets is bounded by

$$\frac{D}{k}\sum_{p=1}^{k}\frac{p}{n}|\psi_p| = O\left(\frac{k}{n}\right),$$

because by Condition C.1, $g(\lambda)$ is continuously differentiable. Next, when $s - [k_1 \log k_1] < q \le s + [k_1 \log k_1]$, the second term of (5.28) is also $o(1)$, since there are at most $O([k_1 \log k_1])$ terms such that $|p + q - s| < [k_1 \log k_1]$, and hence by Lemma 6.1, $D^{-1}\log k_1 \le \log(\tilde{f}_{p+q}\{\lambda_{|p+q-s|_+}^{\alpha} + f_{p+q}^{-1}\}) \le D$, whereas for the remaining ones $|p + q - s| > [k_1 \log k_1]$, so that proceeding as before, it will be $o(1)$ by Condition C.4. Thus, we conclude that the second term of (5.28) in this region is $O(k^{-1}k_1 \log k_1) + o(1) = o(1)$ by Condition C.4. Similarly, when $\max\{0, s - [k_1 \log k_1] - k\} \le q \le s - [k_1 \log k_1]$ we obtain that the second term of (5.28) is also $o(1)$. Proceeding as with the proof of the second term of (5.28), it follows that the third term of (5.28) is $o(1)$ uniformly in $q$.

Using $\int_0^1 \psi(x)\, dx = 0$ and Lemma 6.10, we conclude that

$$\sup_{0\le q\le[n/2]}|\hat{\alpha}(\lambda_q) - \tfrac{1}{2}\delta_{+,n}(\lambda_q) - \tfrac{1}{2}\delta_{-,n}(\lambda_q)| = o_p(1),$$

where $\delta_{+,n}(\lambda_q) = -\alpha h_\psi^{-1} k^{-1}\sum_{p=1}^{k}\psi_p\log(|p+q-s|_+/k)$ and $\delta_{-,n}(\lambda_q) = -\alpha h_\psi^{-1} k^{-1}\sum_{p=1}^{k}\psi_p\log(|q-p-s|_+/k)$.

We now examine the properties of $\delta_{+,n}(\lambda_q)$; those of $\delta_{-,n}(\lambda_q)$ are handled similarly. First, by Lemma 6.10,

$$\delta_{+,n}(\lambda_s) - \alpha = -\alpha\frac{h_\psi^{-1}}{k}\sum_{p=1}^{k}\left(\psi_p\log\left(\frac{p}{k}\right) + h_\psi\right) = O(k^{-1}).$$

Now, for arbitrarily small $\rho > 0$, $\sup_{\rho^{-1}k\le|q-s|}\delta_{+,n}(\lambda_q) < D\rho$ since by Taylor expansion,

$$\sup_{\rho^{-1}k\le|q-s|}|\log(|q-s|_+) - \log(|\pm p + q - s|_+)| \le D\rho.$$

Next, by Proposition 5.1, $\sup_{|q-s|\le\rho k}\delta_{+,n}(\lambda_q) - \alpha < -D\rho^2$, whereas since $\delta_{+,n}(\lambda_q)$ is a nonincreasing function in $|q-s|$, $\sup_{\rho k\le|q-s|\le\rho^{-1}k}\delta_{+,n}(\lambda_q) - \alpha < -D\rho^2$.

Therefore, writing $\widehat{\Lambda}_n(t) = \hat{\alpha}(\lambda_s + (2\pi t)/n)$, we conclude that

$$\Pr\left( \sup_{|t| \geq \rho k} (\widehat{\Lambda}_n(t) - \widehat{\Lambda}_n(0)) > 0 \right) \to 0,$$

that is, $\hat{\lambda}^0$ is a consistent estimator of $\lambda^0$. Thus, to complete the proof of the theorem, it suffices to show that for any $\varepsilon > 0$, there exists $\overline{L} > 0$ such that

(5.30) $$\Pr\left( \sup_{\rho k > |t| > k^{1/2}\overline{L}} (\widehat{\Lambda}_n(t) - \widehat{\Lambda}_n(0)) > 0 \right) < \varepsilon.$$

By Theorem 3.2 (cf. Propositions 5.1–5.3),

$$\frac{k}{t}(\widehat{\Lambda}_n(t) - \widehat{\Lambda}_n(0)) = -\bar{\psi}'' \frac{t\alpha h_\psi^{-1}}{k}\left(1 + O\left(\left(\frac{t}{k}\right)^{1/2}\right)\right) + \frac{1}{t}(\tilde{b}_t + \tilde{\tilde{b}}_t)$$

$$+ \frac{1}{t}(b_{3t} - \tilde{b}_t) + \frac{1}{t}(\hat{\xi}_n^{(5)}(t) - b_{3t})$$

$$+ \frac{1}{t}(b_{4t} - \tilde{\tilde{b}}_t) + \frac{1}{t}(\hat{\xi}_n^{(6)}(t) - b_{4t}).$$

By Proposition 5.3 [cf. (5.19) and (5.20)] and Lemma 6.9,

$$\sup_{|t|<\rho k} \left| \frac{1}{t}(\tilde{b}_t + \tilde{\tilde{b}}_t) \right| = \frac{1}{k^{1/2}} \sup_{|t|<\rho k} \left| \frac{1}{k^{1/2}} \sum_{p=|t|}^{k-|t|} \psi'_p(2\pi I_{\varepsilon,p+s} - 1) + o_p(1) \right|$$

$$= O_p\left(\frac{1}{k^{1/2}}\right),$$

since by Condition C.5 $\psi'(u)$ is continuous, so that $\sup_{|t|<\rho k} |k^{-1/2} \sum_{p=|t|}^{k-|t|} \psi'_p \times (2\pi I_{\varepsilon,p+s} - 1)| = O_p(1)$ by Lemma 6.7, and thence by Lemma 6.5,

(5.31) $$\sup_{|t|<\rho k} \left| \frac{1}{t}(\tilde{b}_t + \tilde{\tilde{b}}_t) + \frac{1}{t}(b_{3t} - \tilde{b}_t) + \frac{1}{t}(b_{4t} - \tilde{\tilde{b}}_t) \right| = O_p(k^{-1/2}).$$

By Condition C.4, there exists a finite positive integer $r$ such that $k_1^{(r-1)\beta} < k^{1/2} < k_1^{(r+1)\beta}$. Consider first the case $2k_1^{1+r\beta} < k$. Then the left-hand side of (5.30) is bounded by

$$\sum_{\ell=1}^{r} \Pr\left( \sup_{L_{\ell-1}k/k_1^{(\ell-1)\beta} \geq |t| > L_\ell k/k_1^{\ell\beta}} \frac{k}{t}(\widehat{\Lambda}_n(t) - \widehat{\Lambda}_n(0)) > 0 \right)$$

(5.32) $$+ \Pr\left( \sup_{L_r k/k_1^{r\beta} \geq |t| > 2k_1 \widetilde{L}} \frac{k}{t}(\widehat{\Lambda}_n(t) - \widehat{\Lambda}_n(0)) > 0 \right)$$

$$+ \Pr\left( \sup_{2k_1 \widetilde{L} \geq |t| > k^{1/2}\overline{L}} \frac{k}{t}(\widehat{\Lambda}_n(t) - \widehat{\Lambda}_n(0)) > 0 \right),$$

where $L_0 = \rho$, $L_\ell > 0$ for $\ell > 1$ and $\widetilde{L} > 0$. Since $h_\psi^{-1}\alpha\bar{\psi}'' > 0$ and $\widehat{\Lambda}_n(t) - \widehat{\Lambda}_n(0) > 0$, the third term of (5.32) is bounded by

$$\Pr\left\{\left|O_p\left(\frac{k_1^2}{k^{3/2}\log k_1} + \frac{k_1}{k_1^\beta k^{1/2}} + 1\right)\right| > \alpha|\bar{\psi}''|\inf_{2k_1\widetilde{L}\geq|t|>k^{1/2}\overline{L}}\left|\frac{t}{k^{1/2}}\right|\right\} < \varepsilon$$

for $\overline{L}$ large enough, by Lemma 6.8(b) and (5.31), and because Condition C.5 implies that $k_1^2 k^{-3/2} = o(1)$, and $k^{-1/2}k_1^{1-\beta} = k^{1/2}k_1^{1+r\beta}/(kk_1^{(1+r)\beta}) = O(1)$. Next, the second term of (5.32) is bounded by

$$\Pr\left\{\left|O_p\left(\frac{k}{k_1^{(r+1)\beta+1}} + \frac{1}{k_1^{r\beta}\log k_1} + \frac{k^{1/2}}{k_1}\right)\right| = o_p(1)$$

$$> \alpha|\bar{\psi}''|\inf_{L_\ell k/k_1^{\ell\beta}>|t|>2\widetilde{L}k_1}\left|\frac{t}{k_1}\right|\right\} < \varepsilon$$

for $\widetilde{L}$ large enough, by Lemma 6.8(a), (5.31) and the fact that $k^{1/2} = o(k_1)$, $k^{1/2} < k_1^{(r+1)\beta}$ and Condition C.4.

Finally, consider the first term of (5.32), whose typical element is, proceeding as before, bounded by

$$\Pr\left\{\left|O_p\left(\frac{k_1^{\ell\beta}}{k_1^{\ell\beta}} + \frac{k_1^{1+\beta}}{k\log k_1} + \frac{k_1^{\ell\beta}}{k^{1/2}}\right)\right| = O_p(1)$$

$$> \alpha|\bar{\psi}''|\inf_{L_{\ell-1}k/k_1^{(\ell-1)\beta}>|t|>L_\ell k/k_1^{\ell\beta}}\left|\frac{tk_1^{\ell\beta}}{k}\right|\right\} < \varepsilon,$$

because $\frac{k_1^{\ell\beta}}{k^{1/2}} \leq \frac{k_1^{\ell\beta+1}}{k^{1/2}k_1} = o(k^{-1}k_1^{\ell\beta+1}) = o(1)$ by Condition C.4 and $\ell \leq r$.

Next consider the case $k < 2k_1^{1+r\beta}$. In this case, let $\tilde{r}$ be the biggest integer such that $2k_1^{1+\tilde{r}\beta} < k < k_1^{1+r\beta}$ but $k < k_1^{1+(\tilde{r}+1)\beta}$. The proof now proceeds identically as in the previous case, but now the sum of the first term of (5.32) runs from $\ell = 1$ to $\tilde{r}$. Note that if $k < 2k_1^{1+\beta}$, so that $r = 1$, then the left-hand side of (5.30) is bounded by

$$\Pr\left(\sup_{\rho k>|t|>2\widetilde{L}k_1}\frac{k}{t}(\widehat{\Lambda}_n(t) - \widehat{\Lambda}_n(0)) > 0\right)$$

$$+ \Pr\left(\sup_{2\widetilde{L}k_1\geq|t|>k^{1/2}\overline{L}}\frac{k}{t}(\widehat{\Lambda}_n(t) - \widehat{\Lambda}_n(0)) > 0\right),$$

and then proceed as in the proof of the second and third terms of (5.32), since Condition C.4 implies that $k_1^3 = o(k^2)$, and recalling that now $k < 2k_1^{1+\beta}$, we have $k_1^{1-\beta}k^{-1/2} = o(1)$ by Condition C.4. □

5.4. *Proof of Theorem* 3.4. We begin with part (a). Observing that Condition C.7 implies that $|\bar{h}_w - h_w| = O(m^{-1})$ and $m^{-1}\sum_{p=1}^{m} w_p = O(m^{-1})$ by Lemma 6.10, the behavior of $(2m)^{1/2}(\check{\alpha}(\lambda_s) - \alpha)$ is governed by that of

$$
\begin{aligned}
&\frac{h_w^{-1}}{(2m)^{1/2}} \sum_{p=1}^{m} w_p \log(\hat{g}_{p+s}\hat{g}_{s-p}) + \frac{h_w^{-1}}{(2m)^{1/2}} \sum_{p=1}^{m} w_p \log\left(\frac{\tilde{f}_{p+s}}{\lambda_p^{-\alpha}} \frac{\tilde{f}_{s-p}}{\lambda_p^{-\alpha}}\right) \\
&\quad - \frac{(2m)^{1/2}\alpha}{h_w}\left(\frac{1}{m}\sum_{p=1}^{m} w_p \log\left(\frac{p}{m}\right) + h_w\right).
\end{aligned}
\tag{5.33}
$$

Recall that $\tilde{f}_p^{-1}\hat{f}_p = \hat{g}_p$. By Lemma 6.10 the last term of (5.33) is $o(1)$. Denoting $m_1^* = [m_1 \log m_1]$, the second term of (5.33) is

$$
\frac{h_w^{-1}}{(2m)^{1/2}}\left\{\sum_{p=1}^{m_1^*} w_p \log\left(\frac{\tilde{f}_{p+s}}{\lambda_p^{-\alpha}} \frac{\tilde{f}_{s-p}}{\lambda_p^{-\alpha}}\right) + \sum_{p=m_1^*+1}^{m} w_p \log\left(\frac{\tilde{f}_{p+s}}{\lambda_p^{-\alpha}} \frac{\tilde{f}_{s-p}}{\lambda_p^{-\alpha}}\right)\right\}.
\tag{5.34}
$$

Next, because by Lemma 6.1 $D^{-1} < \lambda_{m_1}^{\alpha} \tilde{f}_{s\pm p} < D$ for $|p| \leq 2m_1^*$, we have that the first term of (5.34) is $O(m^{-1/2}\log(m_1)\sum_{p=1}^{m_1^*}|w_p|) = O(m_1^{1+\zeta}m^{-(2\zeta+1)/2} \times \log^{1+\zeta}(m_1)) = o(1)$ by Conditions C.6 and C.7 since $\zeta \geq 1/3$. Denoting $g(\lambda_p) = g_p$, the second term of (5.34) is

$$
\begin{aligned}
&\frac{h_w^{-1}}{(2m)^{1/2}} \sum_{p=m_1^*+1}^{m} w_p \log\left(\frac{\tilde{f}_{p+s}}{f_{p+s}} \frac{\tilde{f}_{s-p}}{f_{s-p}}\right) + \frac{h_w^{-1}}{(2m)^{1/2}} \sum_{p=m_1^*+1}^{m} w_p \log(g_{p+s}g_{s-p}) \\
&\quad - \frac{\alpha h_w^{-1}}{(2m)^{1/2}} \sum_{p=m_1^*+1}^{m} w_p\{\log(\lambda_p^{-1}|\lambda_p + \lambda_s - \lambda^0|) + \log(\lambda_p^{-1}|\lambda_s - \lambda^0 - \lambda_p|)\}.
\end{aligned}
$$

Because, for $|p| > 2m_1^*$, Lemma 6.1(a) implies that $D^{-1} \leq m_1^{-2}p^2|f_{p+s}^{-1}\tilde{f}_{p+s} - 1| \leq D$, we obtain by the mean value theorem that $\log f_{p+s}^{-1}\tilde{f}_{p+s} = O(m_1^2 p^{-2})$ and so the first term of the last displayed expression is bounded in absolute value by

$$
\frac{Dm_1^2}{(2m)^{1/2}} \sum_{p=m_1^*+1}^{m} |w_p|p^{-2} = o(m_1^{1+\zeta}m^{-(2\zeta+1)/2}),
$$

whereas the second term is $4\pi^2 Bc^{5/2}/(2^{1/2}h_w) + O(m^{-1/2})$ because

$$
\begin{aligned}
&\frac{1}{(2m)^{1/2}} \sum_{p=m_1^*+1}^{m} w_p \log(g_{p+s}g_{s-p}) \\
&= \frac{2}{(2m)^{1/2}} \sum_{p=1}^{m_1^*} w_p \log(g_s) + \frac{4\pi^2 Bc^{5/2}}{2^{1/2}h_w} + O\left(\frac{1}{m^{1/2}}\right)
\end{aligned}
$$

by Condition C.6 and Taylor expansion of $\log(g_{p+s})$ and $\log(g_{s-p})$ around $\log(g_s)$ and that by Lemma 6.10 and Condition C.7, $\sum_{p=1}^{m} w_p = O(1)$. So, the second term of (5.33) is $4\pi^2 Bc^{5/2}/(2^{1/2}h_w) + o(1)$. Finally, proceeding as with the second term of (5.4), the third term is easily shown to be bounded by $Dm^{-1/2} \sum_p |w_p| p^{-1} = o(1)$.

Denoting $\vartheta_b = \hat{g}_b - 1$, and proceeding as with the proof of Lemma 6.5, to complete the proof of the theorem, it suffices to show that

$$(5.35) \qquad \frac{h_w^{-1}}{(2m)^{1/2}} \sum_{p=1}^{m} w_p(\vartheta_{p+s} + \vartheta_{s-p}) \xrightarrow{d} \mathcal{N}(0, h_w^{-2}\Phi^2),$$

$$(5.36) \qquad \frac{1}{(2m)^{1/2}} \left\{ \sum_{p=1}^{2m_1} + \sum_{p=2m_1+1}^{m} \right\} w_p(\vartheta_{p+s}^2 + \vartheta_{s-p}^2) \xrightarrow{P} 0.$$

We begin with (5.36). By Lemma 6.2(b), the first moment of the first sum inside the braces on the left-hand side of (5.36) is $o(m^{-1/2} \sum_{p=1}^{2m_1} |w_p|) = o(m_1^{1+\zeta} m^{-(2\zeta+1)/2}) = o(1)$ by Conditions C.7 and C.6 since $\zeta \geq 1/3$, whereas the contribution due to the second sum inside the braces on the left-hand side of (5.36) is $O_p(m_1^{-1} m^{1/2}) = o_p(1)$ by Lemma 6.2(a), Condition C.6 and Markov's inequality. So, it remains to show (5.35), whose left-hand side is

$$(5.37) \qquad \frac{h_w^{-1}}{(2m)^{1/2}} \left( \sum_{p=1}^{2m_1} w_p(\vartheta_{p+s} + \vartheta_{s-p}) + \sum_{p=2m_1+1}^{m} w_p(\vartheta_{p+s} + \vartheta_{s-p}) \right).$$

Because Lemma 6.2(b) implies that $E|\vartheta_b| = o(1)$ for $|b - s| < 2m_1$, the first term of (5.37) is $o_p(m_1^{1+\zeta} m^{-(2\zeta+1)/2}) = o_p(1)$ by Conditions C.7 and C.6 and Markov's inequality. But proceeding as in the proof of Proposition 5.3 [cf. (5.21)], the second term of (5.37) converges in distribution to $\mathcal{N}(0, h_w^{-2}\Phi^2)$, which completes the proof of part (a).

Part (b). Dropping the constant $h_w$ and $(2m)^{-1/2}$, it suffices to show that

$$(5.38) \qquad \sum_{p=1}^{m} w_p\{\log(\hat{g}_{p+s-t}/\hat{g}_{p+s}) + \log(\hat{g}_{s-p-t}/\hat{g}_{s-p})\} = o_p(m^{1/2})$$

holds uniformly in $|t| \leq [k^{1/2}M] = o(m^{1/2})$. We only examine the contribution due to the first term on the left-hand side of (5.38); the contribution due to the second term follows by identical steps. The first term on the left-hand side of (5.38) is

$$(5.39) \qquad \sum_{p=1}^{2t} w_p \log(\hat{g}_{p+s-t}/\hat{g}_{p+s}) + \sum_{p=2t+1}^{m} w_p \log(\hat{g}_{p+s-t}/\hat{g}_{p+s}).$$

Using (5.7)–(5.9) and Markov's inequality, the first term of (5.39) is, uniformly in $t$, $o_p(\sum_{p=1}^{[2k^{1/2}M]} |w_p|) = o_p(m^{1/2})$ by Condition C.6. Next, the second term of (5.39) is

$$\sum_{p=t+1}^{2t} w_{p+t} \log \hat{g}_{p+s} - \sum_{p=m-t+1}^{m} w_p \log \hat{g}_{p+s} + \sum_{p=2t+1}^{m-t} (w_{p+t} - w_p) \log \hat{g}_{p+s}.$$

By (5.7)–(5.9), the first two terms of the last displayed expression, uniformly in $t$, are $o_p(k^{1/2}) = o_p(m^{1/2})$ by Condition C.6.

Finally, we consider the third term in the last displayed expression. Let $\vartheta_p = \hat{g}_p - 1$. Since by Lemma 6.3 and Markov's inequality, $\sup_{p=1,\ldots,[n/2]} |\vartheta_p| = o_p(1)$, except in a set $\Omega_n$ such that $\lim_{n \to \infty} \Pr\{\Omega_n\} = 0$, it implies that $\log \hat{g}_{p+s} = \vartheta_{p+s} - 2^{-1} \vartheta_{p+s}^2 (1 + o_p(1))$ by Taylor expansion. So, the third term of the last displayed expression is

(5.40)
$$\left\{ \sum_{p=2m_1+1}^{m-t} + \sum_{p=2t+1}^{2m_1} \right\} (w_{p+t} - w_p) \vartheta_{p+s}$$
$$+ D \sum_{p=2t+1}^{m-t} |w_{p+t} - w_p| \vartheta_{p+s}^2 (1 + o_p(1)).$$

Since Condition C.7 implies that $|w_{p+t} - w_p| \leq D(t/m)^\varsigma$, we have that

$$\sup_{t \leq [k^{1/2}M]} \left| \sum_{p=2t+1}^{m-t} |w_{p+t} - w_p| \vartheta_{p+s}^2 \right|$$
$$\leq D \sup_{t \leq [k^{1/2}M]} \left(\frac{t}{m}\right)^\varsigma \left( \sum_{p=2m_1+1}^{m} \vartheta_{p+s}^2 \right) + \sup_{t \leq [k^{1/2}M]} \left| \sum_{p=2t+1}^{2m_1} w_p \vartheta_{p+s}^2 \right|.$$

But by Lemma 6.2, $E|\vartheta_{p+s}|^2 = O(m_1^{-1})$ if $|p| > 2m_1$, whereas Lemma 6.3 implies that $\sup_{p=1,\ldots,[n/2]} \vartheta_{p+s}^2 = o_p(1)$. Hence, by Markov's inequality and Conditions C.6 and C.7, the third term of (5.40) is $o_p(k^{\varsigma/2} m^{1-\varsigma/2} m_1^{-1} + m^{-\varsigma} m_1^{1+\varsigma}) = o_p(m^{1/2})$. Proceeding similarly and in view of Condition C.7, the second term of (5.40) is $o_p(m_1^{1/2}(k^{1/2}/m)^\varsigma) = o_p(m^{1/2})$. So, denoting $b_t = \sum_{p=2m_1+1}^{m-t} (w_{p+t} - w_p) \vartheta_{p+s}$, that is, the first term of (5.40), to complete the proof we need to show that $\sup_{t \leq [k^{1/2}M]} |b_t| = \sup_{q=1,\ldots,[Mk^{1/4}]} \sup_{(q-1)k^{1/4} \leq t \leq qk^{1/4}} |b_t|$ is $o_p(m^{1/2})$. Now, by the triangle inequality, $\sup_{t \leq [k^{1/2}M]} |b_t|$ is bounded by

$$\sup_{q=1,\ldots,[Mk^{1/4}]} \sup_{(q-1)k^{1/4} \leq t \leq qk^{1/4}} \left| \sum_{p=2m_1+1}^{m-t} (w_{p+t} - w_{p+qk^{1/4}}) \vartheta_{p+s} \right|$$

(5.41) $+ \sup_{q=1,\ldots,[Mk^{1/4}]} \sup_{(q-1)k^{1/4} \leq t \leq qk^{1/4}} \left| \left\{ \sum_{p=2m_1+1}^{m} \right. \right.$

$$- \sum_{p=m-t+1}^{m} \Bigg\} (w_{p+qk^{1/4}} - w_p) \vartheta_{p+s} \Bigg|.$$

Because $(\sup_j |c_j|)^\mu = \sup_j |c_j|^\mu \leq \sum_j |c_j|^\mu$ for $\mu > 0$, the second moment of the second term of (5.41) is bounded by

$$\sum_{q=1}^{[Mk^{1/4}]} \Bigg\{ E \Bigg| \sum_{p=2m_1+1}^{m} (w_{p+qk^{1/4}} - w_p) \vartheta_{p+s} \Bigg|^2$$
$$+ \sum_{\ell=1}^{[Mk^{1/2}]} E \Bigg| \sum_{p=m-\ell+1}^{m} (w_{p+qk^{1/4}} - w_p) \vartheta_{p+s} \Bigg|^2 \Bigg\}$$
$$\leq D \sum_{q=1}^{[Mk^{1/4}]} \Bigg\{ m + \sum_{\ell=1}^{[Mk^{1/2}]} \ell \Bigg\} \bigg(\frac{qk^{1/4}}{m}\bigg)^{2\zeta} = D \frac{k^{\zeta+1/4}}{m^{2\zeta-1}} = o(m),$$

proceeding as with the proof of (5.16) and noting that Condition C.7 implies that $|(m/p)^\zeta w_p| \leq D$ and $|(m/(qk^{1/4}))^\zeta (w_{p+qk^{1/4}} - w_p)| \leq D$ with $\zeta \geq 1/3$. The second moment of the first term of (5.41) is bounded by

$$\sum_{q=1}^{[Mk^{1/4}]} \sum_{t=(q-1)k^{1/4}+1}^{qk^{1/4}} E \Bigg| \sum_{p=2m_1+1}^{m-t} (w_{p+t} - w_{p+qk^{1/4}}) \vartheta_{p+s} \Bigg|^2$$
$$= \frac{k^{\zeta/2}}{m^{2\zeta}} \sum_{q=1}^{[Mk^{1/4}]} \sum_{t=(q-1)k^{1/4}+1}^{qk^{1/4}} E \Bigg| \sum_{p=2m_1+1}^{m-t} w^*_{p,t} \vartheta_{p+s} \Bigg|^2,$$

where $|w^*_{p,t}| = |(m/k^{1/4})^\zeta (w_{p+t} - w_{p+qk^{1/4}})| \leq D$ by Condition C.7. Now proceeding as with the proof of (5.16), the right-hand side of the last displayed equation is $O(m^{1-2\zeta} k^{(1+\zeta)/2}) = o(m)$ since $\zeta \geq 1/3$ and $k = o(m)$ by Condition C.6. Using Markov's inequality we conclude that (5.41) is $o_p(m^{1/2})$ and the proof is complete. $\square$

**6. Technical lemmas.** From now on $\sum_j$ denotes $\sum_{j=-k_1}^{k_1}$ and $k_1 n^{-1} \to 0$.

LEMMA 6.1. *Let $\tilde{f}_p$ be as defined in the proof of Theorem 3.2. Then*

(a)
$$D^{-1} < f_p^{-1} \tilde{f}_p < D,$$
$$(f_p^{-1} \tilde{f}_p - 1) = O(k_1^2/|p-s|^2), \qquad |p-s| \geq 2k_1.$$

(b)
$$D^{-1} \leq \lambda_{k_1}^\alpha \tilde{f}_p \leq D \qquad if\ |p-s| < 2k_1.$$

PROOF. First observe that by definition of $\tilde{f}_p$, $f_p^{-1}\tilde{f}_p = (2k_1+1)^{-1}\sum_j f_p^{-1} \times f_{j+p}$. We begin with (a). We first show that $D^{-1} < (2k_1+1)^{-1}\sum_j f_p^{-1}f_{j+p} < D$. Because $|\lambda^0 - \lambda_s| \le \frac{\pi}{n}$, $D^{-1} < |1 + \frac{n|\lambda^0-\lambda_s|}{2\pi(j+p-s)}| < D$, for $|j| \le k_1$, so that Condition C.1 implies that

$$\frac{D^{-1}}{k_1} \sum_{j=[k_1/4]}^{[k_1/2]} \left|\frac{p-s}{j+p-s}\right|^\alpha \le \frac{1}{2k_1+1}\sum_j \frac{f_{j+p}}{f_p} \le \frac{D}{k_1}\sum_j \left|\frac{p-s}{j+p-s}\right|^\alpha.$$

But $|p-s| \ge 2k_1$ and $|j| \le k_1$ imply that $2/3 < |(p-s)/(j+p-s)| < 2$. From here the conclusion is standard since $\alpha > 0$; we conclude the first part of (a). Next, we show the second part of (a). By Taylor expansion of $f_{j+p}$, the left-hand side is

$$\frac{1}{2k_1+1}\sum_j \left\{\frac{(2\pi)j}{n}\frac{f_p'}{f_p} + \frac{(2\pi)^2 j^2}{2n^2}\frac{f''(\bar{\lambda})}{f_p}\right\} \le \frac{D}{k_1}\sum_j \left(\frac{j}{p-s+\delta j}\right)^2 (1 + o(1)),$$

where $\bar{\lambda} = \bar{\lambda}(p+\delta j)$ is an intermediate point between $\lambda_p$ and $\lambda_{p+j}$ and $\delta = \delta(j) \in (0,1)$, by Condition C.1 and the fact that $f_p^{-1}f(\bar{\lambda})$ is bounded. The conclusion follows since $|p-s+\delta j| \ge |p-s| - |\delta j| > |p-s|/2$.

(b) It is immediate since by Condition C.1, $f_{(j+p)\mathcal{I}(j+p\ne s)} = D\lambda_{|j+p-s|_+}^{-\alpha}(1+o(1))$ and $\lambda_1^{-\alpha}\lambda_{k_1}^\alpha = o(k_1)$. □

LEMMA 6.2. *Denote $\varphi(k_1) = O(k_1^{-1/2}) + O(k_1^{\alpha-1})\mathcal{I}(\alpha \ge 1/2)\mathcal{I}(|p-s| < 2k_1)$. Then*

(a) $$E|\tilde{f}_p^{-1}\ddot{f}_p - 1| = \varphi(k_1),$$

(b) $$E|\tilde{f}_p^{-1}\hat{f}_p - 1| = \varphi(k_1).$$

PROOF. We begin with (a). $\tilde{f}_p^{-1}\ddot{f}_p - 1$ is

$$\tilde{f}_p^{-1}(2k_1+1)^{-1}\sum_{j+p\ne s} f_{j+p}\left(\frac{I_{j+p}}{f_{j+p}} - 1\right)$$
(6.1)
$$+ \tilde{f}_p^{-1}(2k_1+1)^{-1}(I_s - f_{s+1})\mathcal{I}(|p-s| \le k_1).$$

In view of Propositions A.1 and A.2 of [22] and Lemma 6.1, the first term of (6.1) is $\varphi(k_1)$, whereas the second term of (6.1) is also $\varphi(k_1)$ by Lemma 6.1(b) and $E(n^{-\alpha}I_s) < D$.

To show part (b), it suffices to examine $\tilde{f}_p^{-1}(\hat{f}_p - \ddot{f}_p)$, which is by definition

(6.2)
$$\tilde{f}_p^{-1}(n^{-1} - \ddot{f}_p)\mathcal{I}(\ddot{f}_p < n^{-1})$$
$$= ((\tilde{f}_p^{-1}n^{-1} - 1) - (\tilde{f}_p^{-1}\ddot{f}_p - 1))\mathcal{I}(\ddot{f}_p < n^{-1}).$$

By the Cauchy–Schwarz inequality, the second moment of the right-hand side of (6.2) is bounded by

$$2E(\tilde{f}_p^{-1}\ddot{f}_p - 1)^2 + 2(\tilde{f}_p^{-1}n^{-1} - 1)^2 E(\mathcal{I}(\ddot{f}_p < n^{-1})) \leq DE(\tilde{f}_p^{-1}\ddot{f}_p - 1)^2,$$

using the fact that $E(\mathcal{I}(\ddot{f}_p < n^{-1}))$ is

$$\Pr\{\tilde{f}_p^{-1}\ddot{f}_p - 1 < \tilde{f}_p^{-1}n^{-1} - 1\} \leq \Pr\{|\tilde{f}_p^{-1}\ddot{f}_p - 1| > |1 - \tilde{f}_p^{-1}n^{-1}|\},$$

because by Lemma 6.1, $\tilde{f}_p^{-1}n^{-1} - 1 < -D$ for $n$ large enough. Now use part (a) and Markov's inequality to conclude. □

LEMMA 6.3. *Let* $2k_1 < v < u \leq [n/2]$ *and* $p = 0, 1, \ldots, [n/2]$. *Denoting* $\psi(v, u) = O(\frac{\max(u-v, k_1)^{1/\tau}}{k_1^{(2+\tau^2)/2\tau^2}})$ *and* $\varphi(k) = O(\log^{-\mu-1} k)$,

(a)
$$E\left(\sup_{p:\,|p-s|\leq 2k_1} |\tilde{f}_p^{-1}(\ddot{f}_p - \tilde{f}_p)|^\mu\right) = \varphi(k_1),$$

$$E\left(\sup_{p:\,2k_1<|p-s|=v+1,\ldots,u} |\tilde{f}_p^{-1}(\ddot{f}_p - \tilde{f}_p)|\right) = \psi(v, u),$$

(b)
$$E \sup_{p:\,|p-s|\leq 2k_1} |\tilde{f}_p^{-1}(\ddot{f}_p - \hat{f}_p)|^\mu = \varphi(k_1),$$

$$E \sup_{p:\,2k_1<|p-s|=v+1,\ldots,u} |\tilde{f}_p^{-1}(\ddot{f}_p - \hat{f}_p)| = \psi(v, u).$$

PROOF. For notational simplicity we shall take $s = 0$. We begin with part (a). From Hidalgo and Robinson's [22] Proposition A.1, it suffices to examine the behavior of $\tilde{f}_p^{-1}(\ddot{f}_p - E\ddot{f}_p)$. On the other hand, Hidalgo and Robinson's [22] Proposition A.3(a) and (b) implies that it suffices to examine the behavior of $\tilde{f}_p^{-1}(\ddot{f}_{\varepsilon,p} - E\ddot{f}_{\varepsilon,p})$, where

$$\ddot{f}_{\varepsilon,p} = \frac{1}{2k_1 + 1} \sum_j f_{j+p} I_{\varepsilon,j+p}$$

and $I_{\varepsilon,p} = I_\varepsilon(\lambda_p)$ is the periodogram of $\{\varepsilon_t\}_{t=1}^n$. We examine $\sup_{p=v+1,\ldots,u} |\tilde{f}_p^{-1} \times (\ddot{f}_{\varepsilon,p} - E\ddot{f}_{\varepsilon,p})|$ only; that of $\sup_{p=1,\ldots,2k_1} |\tilde{f}_p^{-1}(\ddot{f}_{\varepsilon,p} - E\ddot{f}_{\varepsilon,p})|$ is similarly handled. Because $\sup_j |a_j| = (\sup_j |a_j|^\tau)^{1/\tau}$, the $\tau$th power of $\sup_{p=v+1,\ldots,u} |\tilde{f}_p^{-1} \times (\ddot{f}_{\varepsilon,p} - E\ddot{f}_{\varepsilon,p})|$ is, except for constants,

$$\sup_{p=v+1,\ldots,u} \left|\frac{1}{2k_1 + 1} \sum_j \phi_{j+p,p}((2\pi)I_{\varepsilon,j+p} - 1)\right|^\tau,$$

where $\phi_{j,p} = \tilde{f}_p^{-1} f_j$. The last displayed expression is bounded by

$$
(6.3) \quad \begin{aligned}
& 2^{\tau-1} \sup_q \sup_p \left| \frac{1}{2k_1+1} \sum_j (\phi_{j+p,p}((2\pi)I_{\varepsilon,j+p} - 1) \right. \\
& \left. \qquad\qquad - \phi_{j+b,p}((2\pi)I_{\varepsilon,j+b} - 1)) \right|^\tau \\
& + 2^{\tau-1} \sup_q \sup_p \left| \frac{1}{2k_1+1} \sum_j \phi_{j+b,p}((2\pi)I_{\varepsilon,j+b} - 1) \right|^\tau,
\end{aligned}
$$

where $\sup_q$ and $\sup_p$ denote $\sup_{q=1+v/k_1^{1/\tau},\ldots,u/k_1^{1/\tau}}$ and $\sup_{p=1+b-k_1^{1/\tau},\ldots,b}$, respectively, and $b = qk_1^{1/\tau}$.

After the change of indices $j = j' - k_1$, the second term of (6.3) is bounded by

$$
(6.4) \quad \begin{aligned}
& D \sup_q \sup_p \left| \frac{1}{2k_1+1} \sum_{j=0}^{2k_1-1} (\phi_{j+b-k_1,p} - \phi_{j+b+1-k_1,p}) \right. \\
& \left. \qquad\qquad \times \sum_{a=0}^{j} ((2\pi)I_{\varepsilon,a+b-k_1} - 1) \right|^\tau \\
& + D \sup_q \sup_p |\phi_{b+k_1,p}|^\tau \left| \frac{1}{2k_1+1} \sum_{j=0}^{2k_1} ((2\pi)I_{\varepsilon,j+b-k_1} - 1) \right|^\tau,
\end{aligned}
$$

by Abel summation by parts. On the other hand, by Conditions C.1 and C.3,

$$
\begin{aligned}
|\phi_{j+b-k_1,p} - \phi_{j+b+1-k_1,p}| &\leq D \tilde{f}_p^{-1} (j+b-k_1)^{-1-\alpha} n^\alpha \\
&\leq D \left( \frac{p}{j+b-k_1} \right)^\alpha (j+b-k_1)^{-1},
\end{aligned}
$$

since by Lemma 6.1(a) $D^{-1} < |f_p^{-1} \tilde{f}_p| < D$. So, using $\sup_j |a_j|^\tau \leq \sum_j |a_j|^\tau$, by Hölder's inequality and $D^{-1} < |\phi_{b+k_1,p}| < D$ by Lemma 6.1(a), we obtain that the first moment of (6.4) is bounded by

$$
\frac{D}{2k_1+1} \sum_{q=1+v/k_1^{1/\tau}}^{u/k_1^{1/\tau}} \sum_{j=0}^{2k_1} \sup_p \left( \frac{p}{j+b-k_1} \right)^{\tau\alpha} (j+b-k_1)^{-\tau}
$$

$$
\times E \left| \sum_{a=0}^{j} ((2\pi)I_{\varepsilon,a+b-k_1} - 1) \right|^\tau
$$

$$+ D \sum_{q=1+v/k_1^{1/\tau}}^{u/k_1^{1/\tau}} E\left|\frac{1}{2k_1+1}\sum_{j=0}^{2k_1}((2\pi)I_{\varepsilon,j+b-k_1}-1)\right|^\tau,$$

which, because $(2\pi)EI_{\varepsilon,j+p}=1$ and proceeding as in the proof of Brillinger's [5] Theorem 7.4.4, is bounded by

$$D \sum_{q=1+v/k_1^{1/\tau}}^{u/k_1^{1/\tau}} \left(\left(\frac{b}{b-k_1}\right)^{\tau\alpha} \frac{1}{2k_1} \sum_{j=0}^{2k_1} \frac{(j+1)^{\tau/2}}{(j+b-k_1)^\tau} + k_1^{-\tau/2}\right)$$

$$\leq D \sum_{q=1+v/k_1^{1/\tau}}^{u/k_1^{1/\tau}} \left(\left(\frac{1}{b}\right)^{\tau/2} + \left(\frac{1}{k_1}\right)^{\tau/2}\right) = O\left(\frac{\max(u-v,k_1)}{k_1^{\tau/2+1/\tau}}\right),$$

because $\alpha < 1$, $b \leq 2(b-k_1)$, $q \geq 1+v/k_1^{1/\tau}$ and $b=qk_1^{1/\tau}$. Thus, we conclude that the second term of (6.3) is $O(\max(u-v,k_1)/k_1^{\tau/2+1/\tau})$.

Next, we examine the first term of (6.3). Because

$$a_{p,k_1} = \sum_{j=0}^{2k_1}(\phi_{j+p-k_1,p}((2\pi)I_{\varepsilon,j+p-k_1}-1) - \phi_{j+b-k_1,p}((2\pi)I_{\varepsilon,j+b-k_1}-1))$$

has at most $k_1^{1/\tau}$ terms, and because $(2\pi)EI_{\varepsilon,j+p}=1$ and proceeding as in the proof of Brillinger's [5] Theorem 7.4.4, its $\tau$th moment is bounded by $k_1^{1/2}$, so that the expectation of the first term of (6.3) is bounded by

$$D \sum_{q=1+v/k_1^{1/\tau}}^{u/k_1^{1/\tau}} \sum_{p=1+b-k_1^{1/\tau}}^{b} k_1^{1/2-\tau} = o_p\left(\frac{u-v}{k_1^{\tau/2+1/\tau}}\right),$$

because $\tau > 2$. This completes the proof of part (a).

To show part (b), denoting $a_p = \tilde{f}_p^{-1}n^{-1} - 1$ and using (6.2), it suffices to examine

$$\sup_p |a_p|\mathcal{I}(\ddot{f}_p < n^{-1}) \leq D\sup_p \mathcal{I}(\tilde{f}_p^{-1}\ddot{f}_p - 1 < a_p).$$

But the expectation of the right-hand side is bounded by

$$DE\mathcal{I}\left(\sup_p |\tilde{f}_p^{-1}\ddot{f}_p - 1| > \min_p |a_p|\right)$$

$$= D\Pr\left\{\sup_p |\tilde{f}_p^{-1}\ddot{f}_p - 1| > \min_p |a_p|\right\}$$

$$\leq D\left(\min_p |a_p|\right)^{-\tau} E\left(\sup_p |\tilde{f}_p^{-1}\ddot{f}_p - 1|\right)^\tau$$

$$\leq D(\varphi(k_1)\mathcal{I}(p \leq 2k_1) + \psi(v,u)\mathcal{I}(p > 2k_1))$$

by Markov's inequality and because $(\sup_p |c_p|)^\tau = \sup_p |c_p|^\tau$. □

LEMMA 6.4. *Let $h(u)$ be a twice continuously differentiable function in $(0,1)$ such that $h(0) = h'(0) = h(1) = 0$, where $h'(u) = \frac{d}{du}h(u)$. Consider a sequence $\{\nu_j\}$ such that $|j^2 \nu_j| \le D$ for all $j$. Then, for $a > 2t$ and denoting $q = p - t + 1$,*

$$(6.5) \quad \sum_{j=a+1}^{p} h(j/p)(\nu_j - \nu_{j-t}) = O\left(\frac{t}{p^2}\log\left(\frac{p}{t}\right) + \frac{t^2}{pq^2} + \frac{t}{p^2}\right).$$

PROOF. The left-hand side of (6.5) is

$$(6.6) \quad \sum_{j=p-t+1}^{p} h(j/p)\nu_j - \sum_{j=a-t+1}^{a} h((j+t)/p)\nu_j + \sum_{j=a+1}^{p-t} (h(j/p) - h((j+t)/p))\nu_j.$$

Since the first derivative of $h(u)$ is continuous and $h(1) = 0$, from the mean value theorem it follows that the absolute value of the first term of (6.6) is bounded by

$$D \sum_{j=p-t+1}^{p} |(p-j)/p||\nu_j| = O(t^2/(p(p-t+1)^2)) = O(p^{-1}q^{-2}t^2),$$

using $|j^2 \nu_j| \le D$. The absolute value of the second term of (6.6) is bounded by

$$D \sum_{j=a-t+1}^{a} \left(\frac{j+t}{p}\right)^2 |\nu_j| = O\left(\frac{t}{p^2}\right),$$

since $h(0) = h'(0) = 0$ and $|j^2 \nu_j| \le D$, whereas the absolute value of the third term of (6.6) is bounded by

$$D\frac{t}{p} \sum_{j=a+1}^{p-t} \left|h'\left(\frac{j}{p}\right) + \frac{t}{p}h''\left(\frac{j}{p} + \xi\frac{t}{p}\right)\right||\nu_j| = O\left(\frac{t}{p^2}\log\left(\frac{p}{t}\right)\right)$$

by Taylor expansion of $h'(x)$ and using $h'(0) = 0$, where $\xi = \xi(j) \in (0,1)$. □

LEMMA 6.5. *Let $\tilde{b}_t$, $\tilde{\tilde{b}}_t$, $b_{3t}$ and $b_{4t}$ be given in (5.16), (5.10) and (5.15), respectively. Then, for $\rho < 1/3$,*

(a) $\sup_{t \le \rho k} t^{-1}|b_{3t} - \tilde{b}_t| = o_p(k^{-1/2})$,

(b) $\sup_{t\leq \rho k} t^{-1}|b_{4t} - \tilde{\tilde{b}}_t| = o_p(k^{-1/2})$.

PROOF. We only examine part (a); part (b) is identical. Because by Lemma 6.3, $\sup_{\ell=1,\ldots,[n/2]} |\hat{g}_\ell - 1| = O_p(\log^{-2} k_1)$, then except in a set $\Omega_n$ such that $\lim_n \Pr\{\Omega_n\} = 0$, $\log \hat{g}_{.+s} = \vartheta_{.+s} - 2^{-1}\vartheta^2_{.+s}(1 + o_p(1))$ by Taylor expansion, which implies that for $n$ sufficiently large, by definition of $b_{3t}$,

$$\sup_{t\leq \rho k} t^{-1}|b_{3t} - \tilde{b}_t| \leq D \sup_{t\leq \rho k} t^{-1} \sum_{p=2t+1}^{k-t} |\psi_p - \psi_{p+t}|\vartheta^2_{p+s}.$$

Because by Lemma 6.2, for $|p| > 2k_1$, $\vartheta^2_{p+s} = O_p(k_1^{-1})$, for $|p| < 2k_1$, $\vartheta^2_{p+s} = O_p(k_1^{2(\alpha-1)}\mathcal{I}(\alpha \geq 1/2) + k_1^{-1})$ and by Condition C.5, $|\psi_p - \psi_{p+t}| \leq D|\psi'_\xi|t/k$, where $\psi'_p = \psi'(p/k)$ and $p \leq \xi \leq p+t$, the last displayed expression is

$$O_p\left(\frac{\mathcal{I}(\alpha \geq 1/2)}{k_1^{2(1-\alpha)}} \sum_{p=1}^{2k_1} \frac{\sup_{t<2k_1} |\psi'_\xi|}{k} + \frac{1}{k_1}\sum_{p=1}^{k} \frac{\sup_{t\leq \rho k} |\psi'_\xi|}{k}\right) = o_p(k^{-1/2})$$

by Conditions C.5 and C.4. □

LEMMA 6.6. *Let $\phi_p = \phi(p/k)$, where $\phi(u)$ is a continuous function in $(0,1)$. Define*

$$c_r(\mu; \vartheta) = \frac{2}{nk^{1/2}} \sum_{p=[k\mu]+1}^{[k\vartheta]} \phi_p \cos(r\lambda_p),$$

*where $0 \leq \mu < \vartheta \leq 1$. For any $\mu < \vartheta_1 < \vartheta_2 \leq 1$, if $k/n \to 0$, then*

(6.7) $$\sum_{r_1=1}^{n-1}\sum_{r_2=1}^{n-r_1} c_{r_2}(\mu; \vartheta_1) c_{r_2}(\mu; \vartheta_2) = \int_\mu^{\vartheta_1} \phi^2(u)\, du\, (1 + o(1)).$$

PROOF. The left-hand side of (6.7) is

(6.8)
$$\frac{4}{n^2 k} \sum_{p_1=[k\mu]+1}^{[k\vartheta_1]} \phi_{p_1} \sum_{p_2=[k\mu]+1}^{[k\vartheta_2]} \phi_{p_2} \sum_{r_1=1}^{n-1}\sum_{r_2=1}^{n-r_1} \cos(r_2 \lambda_{p_1})\cos(r_2 \lambda_{p_2})$$
$$= \frac{4}{n^2 k} \sum_{p=[k\mu]+1}^{[k\vartheta_1]} \phi_p^2 \sum_{r_1=1}^{n-1}\sum_{r_2=1}^{n-r_1} \cos^2(r_2 \lambda_p)$$
$$+ \frac{2}{n^2 k} \sum_{p_1=[k\mu]+1}^{[k\vartheta_1]} \phi_{p_1} \sum_{p_2=[k\mu]+1, p_2\neq p_1}^{[k\vartheta_2]} \phi_{p_2} \sum_{r_1=1}^{n-1}\sum_{r_2=1}^{n-r_1} \{\cos(r_2 \lambda_{p_1+p_2})$$
$$+ \cos(r_2\lambda_{p_1-p_2})\}.$$

Because (see [32]) $\sum_{r_1=1}^{n-1}\sum_{r_2=1}^{n-r_1}\cos^2(r_2\lambda_p) = (n-1)^2/4$ and

$$\sum_{r_1=1}^{n-1}\sum_{r_2=1}^{n-r_1}\{\cos(r_2\lambda_{p_1+p_2}) + \cos(r_2\lambda_{p_1-p_2})\} = -n \qquad \text{for } p_1 \neq p_2,$$

the right-hand side of (6.8) is

$$\frac{(n-1)^2}{n^2}\left(\frac{1}{k}\sum_{p=[k\mu]+1}^{[k\vartheta_1]}\phi_p^2\right) - \frac{2}{nk}\sum_{p_1=[k\mu]+1}^{[k\vartheta_1]}\phi_{p_1}\sum_{p_2=[k\mu]+1, p_2\neq p_1}^{[k\vartheta_2]}\phi_{p_2}$$

$$= \int_\mu^{\vartheta_1}\phi^2(u)\,du\,(1+o(1)),$$

because $\phi(u)$ is continuous in $u$ and $k/n \to 0$. □

LEMMA 6.7. *Denote $\eta_p = (2\pi)I_{\varepsilon,p} - 1$ and $\phi(u)$ as in Lemma 6.6. The process*

$$R_n(\vartheta) = \frac{1}{k^{1/2}}\sum_{p=[k\vartheta]+1}^{k-[k\vartheta]}\phi_p\eta_p, \qquad 0 \leq \vartheta \leq 1/2,$$

*is tight.*

PROOF. Since by Proposition 5.3, the finite limit distributions of $R_n(\vartheta)$ converge to those of a Gaussian process with continuous paths, then by Billingsley's [4] Theorem 15.6, it suffices to check the moment condition

(6.9) $\qquad E(|R_n(\vartheta_2) - R_n(\vartheta)|^\tau |R_n(\vartheta) - R_n(\vartheta_1)|^\tau) \leq D(\vartheta_2 - \vartheta_1)^\psi$

for some $\tau > 0, \psi > 1$, where $0 \leq \vartheta_1 < \vartheta < \vartheta_2 \leq 1/2$. Because

$$R_n(\vartheta) - R_n(\vartheta_2) = \frac{1}{k^{1/2}}\sum_{p=[k\vartheta]+1}^{[k\vartheta_2]}\phi_p\eta_p + \frac{1}{k^{1/2}}\sum_{p=k-[k\vartheta_2]+1}^{k-[k\vartheta]}\phi_p\eta_p,$$

a sufficient condition for (6.9) to hold is

(6.10)
$$E\left(\left|\frac{1}{k^{1/2}}\sum_{p=[k\vartheta]+1}^{[k\vartheta_2]}\phi_p\eta_p\right|^\tau\left|\frac{1}{k^{1/2}}\sum_{p=[k\vartheta_1]+1}^{[k\vartheta]}\phi_p\eta_p\right|^\tau\right) \leq D(\vartheta_2-\vartheta_1)^\psi,$$

$$E\left(\left|\frac{1}{k^{1/2}}\sum_{p=k-[k\vartheta_2]+1}^{k-[k\vartheta]}\phi_p\eta_p\right|^\tau\left|\frac{1}{k^{1/2}}\sum_{p=k-[k\vartheta]+1}^{k-[k\vartheta_1]}\phi_p\eta_p\right|^\tau\right) \leq D(\vartheta_2-\vartheta_1)^\psi.$$

We will examine the first inequality of (6.10) only; the second displayed inequality is similarly handled.

By definition of $\eta_p$,

$$\frac{1}{k^{1/2}} \sum_{p=[k\vartheta]+1}^{[k\vartheta_2]} \phi_p \eta_p = \left(\frac{1}{k} \sum_{p=[k\vartheta]+1}^{[k\vartheta_2]} \phi_p\right)\left(\frac{k^{1/2}}{n} \sum_{r=1}^{n}(\varepsilon_r^2 - 1)\right)$$

$$+ \sum_{r=2}^{n} \varepsilon_r \sum_{a=1}^{r-1} \varepsilon_a c_{r-a}(\vartheta, \vartheta_2)$$

$$:= \mathcal{E}_{1,n}(\vartheta, \vartheta_2) + \mathcal{E}_{2,n}(\vartheta, \vartheta_2),$$

where $c_r(\vartheta, \vartheta_2)$ was defined in Lemma 6.6. Because $|\sum_{p=[k\vartheta]+1}^{[k\vartheta_2]} \phi_p| \leq Dk|\vartheta_2 - \vartheta|$ by continuity of $\phi(x)$ and $E(\sum_{r=1}^{n}(\varepsilon_r^2 - 1))^2 < Dn$ by Condition C.2, $E(|\mathcal{E}_{1,n}(\vartheta, \vartheta_2)||\mathcal{E}_{1,n}(\vartheta_1, \vartheta)|) \leq (\vartheta_2 - \vartheta_1)^2$ by the Cauchy–Schwarz inequality and $|\vartheta_2 - \vartheta||\vartheta - \vartheta_1| < |\vartheta_2 - \vartheta_1|^2$. That is, $\mathcal{E}_{1,n}(\vartheta, \vartheta_2)$ satisfies the first inequality in (6.10) with $\tau = 1$ and $\psi = 2$. So, to complete the proof, it suffices to examine that the first inequality in (6.10) holds for $\mathcal{E}_{2,n}(\vartheta, \vartheta_2)$. The fourth moment of $\mathcal{E}_{2,n}(\vartheta, \vartheta_2)$ is

$$E\left[\sum_{2=r_1 \leq r_2 \leq r_3 \leq r_4} \prod_{j=1}^{4} \varepsilon_{r_j}\left(\sum_{a_j=1}^{r_j-1} \varepsilon_{a_j} c_{r_j-a_j}(\vartheta, \vartheta_2)\right)\right]$$

$$\leq D \prod_{j=1}^{4} \left(\sum_{1 \leq a_j \leq r_j \leq n} c_{r_j-a_j}^2(\vartheta, \vartheta_2)\right)^{1/2}$$

$$= D \left(\sum_{1 \leq a \leq r \leq n} c_{r-a}^2(\vartheta, \vartheta_2)\right)^2,$$

proceeding as in the proof of Lemma 5.4 of [12]. But proceeding as in Lemma 6.6, the right-hand side of the last displayed equation is bounded by

$$D\left(\int_{\vartheta}^{\vartheta_2} \phi^2(u)\, du\right)^2 \leq D(\vartheta_2 - \vartheta_1)^2$$

since $\phi(u)$ is continuous. This concludes the proof, choosing $\tau = \psi = 2$. □

LEMMA 6.8. *Let $2k_1 \leq t_0 < \rho k$ for some arbitrarily small $\rho > 0$. Then*

(a)

$$\sup_{2k_1 < t \leq t_0} |t^{-1}(\hat{\xi}_n^{(5)}(t) - b_{3t})| = O_p\left(\frac{t_0}{kk_1^\beta} + \frac{k_1 t_0}{k^2 \log k_1}\right),$$

(6.11)

$$\sup_{2k_1 < t \leq t_0} |t^{-1}(\hat{\xi}_n^{(6)}(t) - b_{4t})| = O_p\left(\frac{t_0}{kk_1^\beta} + \frac{k_1 t_0}{k^2 \log k_1}\right),$$

(b)

$$\sup_{Lk^{1/2}<t\leq 2k_1} |t^{-1}(\hat{\xi}_n^{(5)}(t) - b_{3t})| = O_p\left(\frac{k_1^2}{k^2 \log k_1} + \frac{k_1}{kk_1^\beta}\right),$$

(6.12)

$$\sup_{Lk^{1/2}<t\leq 2k_1} |t^{-1}(\hat{\xi}_n^{(6)}(t) - b_{4t})| = O_p\left(\frac{k_1^2}{k^2 \log k_1} + \frac{k_1}{kk_1^\beta}\right),$$

where $b_{3t}$ and $b_{4t}$ are given by (5.10) and (5.15), respectively, and $L > 0$ and $\beta > 0$.

PROOF. We begin with (a). We only examine the first equality in (6.11); the proof of the second equality is similarly handled. By definition (see Proposition 5.3) $\hat{\xi}_n^{(5)}(t) = b_{1t} + b_{2t}$. First, $\sup_{2k_1<t\leq t_0} |t^{-1}b_{1t}|$ satisfies the equality in (6.11) by Condition C.5 and since the sum in $p$ has at most $2k_1$ terms, say $p^* = 1, \ldots, 2k_1$, for which $\sup_{p^*=1,\ldots,2k_1} |\log \hat{g}_{p^*+s}| = O_p(\log^{-1} k_1)$, whereas for the remaining terms $\sup_{p=1,\ldots,k;p\neq p^*} |\log \hat{g}_{p+s}| = O_p(k_1^{-\beta})$ using (5.7) and (5.9).

Next we estimate $b_{2t} - b_{3t}$. First by (5.11), $\sup_{2k_1<t\leq t_0} t^{-1}|b_{2t} - b_{3t}|$ is

$$\sup_{2k_1<t\leq t_0} t^{-1} \left| \sum_{p=k-t+1}^{k} \psi_p \log(\hat{g}_{p+s}) - \sum_{p=t+1}^{2t} \psi_{p+t} \log(\hat{g}_{p+s}) \right|$$

$$= O_p(t_0 k^{-1} k_1^{-\beta} + t_0^2 k^{-2} k_1^{-\beta}),$$

by Condition C.5 and using (5.7) and (5.9) and Markov's inequality.

Part (b). As was done in part (a), we only examine the first equality in (6.12); the second is similarly handled. By Condition C.5 and (5.7) and (5.9) together with Markov's inequality, it follows easily that $\sup_{Lk^{1/2}<t\leq 2k_1} |t^{-1}b_{1t}| = O_p(k^{-2}k_1^2 \log^{-1} k_1)$. Finally, we estimate $b_{2t} - b_{3t}$. As was done in part (a), it suffices to examine

$$\sup_{Lk^{1/2}<t\leq 2k_1} t^{-1} \left\{ \left| \sum_{p=k-t+1}^{k} \psi_p \log(\hat{g}_{p+s}) \right| + \left| \sum_{p=t+1}^{2t} \psi_{p+t} \log(\hat{g}_{p+s}) \right| \right\}.$$

The second term of the last displayed expression is $O_p(k_1^2 k^{-2} \log^{-1} k_1)$ by Condition C.5 and using (5.7)–(5.8), whereas the first term is $O_p(k_1 k^{-1} k_1^{-\beta})$ using Condition C.5 and (5.7) and (5.9), which concludes the proof. □

LEMMA 6.9. *Let $\phi(u)$ be as in Lemma 6.6. Then*

$$\sup_{\omega\in[0,1]} \left| \frac{1}{k^{1/2}} \sum_{j=1}^{[k\omega]} \phi_j \left( \frac{I_{j+s}}{f_{j+s}} - 2\pi I_{\varepsilon,j+s} \right) \right| = o_p(1).$$

PROOF. Writing $u_j = f_{j+s}^{-1/2} \omega_{j+s,x}$ and $v_j = (2\pi)^{1/2} \omega_{j+s,\varepsilon}$ where $\omega_{j+s,x}$ and $\omega_{j+s,\varepsilon}$ are the discrete Fourier transforms of $x_r$ and $\varepsilon_r$, respectively, the left-hand side of the last displayed expression is, by the triangle inequality, bounded by

$$(6.13) \quad \sup_{\omega \in [0,1]} \frac{1}{k^{1/2}} \sum_{j=1}^{[k\omega]} |\phi_j| |u_j - v_j|^2 + 2 \sup_{\omega \in [0,1]} \left| \frac{1}{k^{1/2}} \sum_{j=1}^{[k\omega]} \phi_j v_j (\bar{u}_j - \bar{v}_j) \right|,$$

where $\bar{c}$ denotes the conjugate of the complex number $c$.

The first term of (6.13) is $o_p(1)$ since its expectation is bounded by

$$k^{-1/2} \sum_{j=1}^{k} |\phi_j| \{ (E|u_j|^2 - 1)$$
$$(6.14) \qquad - (E(u_j \bar{v}_j) - 1) - (E(\bar{u}_j v_j) - 1) + (E|v_j|^2 - 1) \}$$
$$= O\left( k^{-1/2} \sum_{j=1}^{k} \frac{\log j}{j} \right),$$

because $E|v_j|^2 = 1$, $|\phi_j| \leq D$ and by the extension of Theorems 1 and 2 of [31] given in Lemma 4.4 of [12].

Next, to show that the second term of (6.13) is $o_p(1)$, it suffices to show that the finite-dimensional distributions of the term inside the absolute value converge to zero and the tightness condition. First, choosing $\omega_1^*$ such that $[k\omega_1^*] = \max([k^\zeta], [k\omega_1])$ for some $0 < \zeta < 1/4$, then for any $0 < \omega_1 < \omega_2 < 1$, $E|k^{-1/2} \sum_{j=1+[k\omega_1]}^{[k\omega_2]} \phi_j v_j (\bar{u}_j - \bar{v}_j)|^2$ is bounded by

$$2E \left| k^{-1/2} \sum_{j=1+[k\omega_1^*]}^{[k\omega_2]} \phi_j v_j (\bar{u}_j - \bar{v}_j) \right|^2 + 2k^{-1} [k\omega_1^*]([k\omega_1^*] - [k\omega_1])$$
$$(6.15) \quad \leq Dk^{-1} \log^2 k ( ([k\omega_2]^{1/2} - [k\omega_1^*]^{1/2})$$
$$\times (\log(k\omega_2) - \log(k\omega_1^*)) + [k\omega_1^*]([k\omega_1^*] - [k\omega_1])),$$

proceeding as with (4.8) in [32]. So, the finite-dimensional distributions of the second term of (6.13) converge to zero in probability by Markov's inequality.

To complete the proof we need to show tightness. Since the limiting process has continuous paths, by Billingsley's [4] Theorem 15.6, it suffices to show that

$$(6.16) \quad E \left| k^{-1/2} \sum_{j=1+[k\omega_1]}^{[k\omega_2]} \phi_j v_j (\bar{u}_j - \bar{v}_j) \right|^4 \leq D(H(\omega_2) - H(\omega_1))^{1+\delta},$$

where $\delta > 0$, $0 < \omega_1 < \omega_2 < 1$ and $H(\omega)$ is a nondecreasing continuous function. The left-hand side of (6.16) is bounded by

$$k^{-2}(|M_4| + 3M_2^2),$$

where $M_r$ denotes the $r$th cumulant of $\sum_{j=1+[k\omega_1]}^{[k\omega_2]} \phi_j v_j(\bar{u}_j - \bar{v}_j)$. Using the inequality in (6.15), $k^{-2}M_2^2 \leq D(H(\omega_2) - H(\omega_1))^{1+\delta}$, so it remains to show that $k^{-2}|M_4|$ satisfies the inequality in (6.16). Now $k^{-2}|M_4|$ is

$$(6.17) \quad \frac{1}{k^2} \sum_{j_1,j_2,j_3,j_4=1+[k\omega_1]}^{[k\omega_2]} \left( \prod_{i=1}^4 \phi_{j_i} \right) \operatorname{cum}(v_{j_1}\bar{z}_{j_1}, v_{j_2}\bar{z}_{j_2}, v_{j_3}\bar{z}_{j_3}, v_{j_4}\bar{z}_{j_4}),$$

where we have abbreviated $u_j - v_j$ by $z_j$. By Theorem 2.3.2 of [5] and denoting $X_{j1} = \phi_j v_j$ and $X_{j2} = \phi_j z_j$,

$$(6.17) = \frac{1}{k^2} \sum_{\vartheta} \operatorname{cum}(X_{j\ell}; j\ell \in \vartheta_1) \cdots \operatorname{cum}(X_{j\ell}; j\ell \in \vartheta_p),$$

where the summation is over all indecomposable partitions $\vartheta = \vartheta_1 \cup \cdots \cup \vartheta_p$. A typical component in $\operatorname{cum}(X_{j\ell}; j\ell \in \vartheta_1)$ has $q_1$ elements $v_j$ and $q_2$ elements $z_j$, so applying formulae of [5], (2.6.3), page 26, and (2.10.3), page 39, we deduce after straightforward calculations that $\operatorname{cum}(X_{j\ell}; j\ell \in \vartheta_1)$ is $\prod_{j \in v_1} \phi_j$ times

$$\frac{\mu_{q_1+q_2}}{k^{(q_1+q_2)/2}}$$

$$\times \int_{[-\pi,\pi]^{q_1+q_2-1}} \frac{\beta(\lambda^1 + \cdots + \lambda^{(q_1-1)} + \nu^1 + \cdots + \nu^{q_2})\beta(-\lambda^1)\cdots\beta(-\lambda^{q_1-1})}{\beta_{j_1}\cdots\beta_{j_{q_1}}}$$

$$\times \tilde{\beta}(-\nu^1)\cdots\tilde{\beta}(-\nu^{q_2})$$

$$\times E_{j_1\cdots j_{q_1}\ell_1\cdots\ell_{q_2}}(\lambda^1,\ldots,\lambda^{(q_1-1)},\nu^1,\ldots,\nu^{q_2})\,d\lambda^1\cdots d\lambda^{(q_1-1)}\,d\nu^1\cdots d\nu^{q_2},$$

where $E_{j_1\cdots j_q\ell_1\cdots\ell_p}(\lambda^1,\ldots,\lambda^{(q-1)},\nu^1,\ldots,\nu^p)$ is

$$G(\lambda_{j_1} - [\lambda^1 + \cdots + \lambda^{(q-1)} + \nu^1 + \cdots + \nu^p])G(\lambda_{j_2} + \lambda^1)$$

$$\times \cdots \times G(\lambda_{j_q} + \lambda^{(q-1)})G(\nu^1 - \lambda_{\ell_1}) \times \cdots \times G(\nu^p - \lambda_{\ell_p}),$$

with $G(\lambda) = \sum_{t=1}^n e^{it\lambda}$ and, say, $\tilde{\beta}(-\nu^1) = \beta_{\ell_1}^{-1}\beta(-\nu^1) - 1$. But by a routine extension of Lemma 3 of [32] and observing that in each partitioned $\upsilon$, the subindex $j_i, i = 1, \ldots, 4$, appears only once,

$$(6.17) \leq Dk^{-2}\left(\sum_{j=1+[k\omega_1]}^{[k\omega_2]} \frac{1}{j^{1/2}}\right)^4 \leq D(H(\omega_2) - H(\omega_1))^4,$$

where $H(\omega) = \omega^{1/2}$, which is a nondecreasing continuous function. $\square$

Remark 6.1. An alternative proof of this lemma can be found in Lemma 4 of [9].

Lemma 6.10 ([5], page 15). *Let $h(x)$, $0 \leq x \leq 1$, be integrable and have an integrable derivative $h^{(1)}(x)$. Then*

$$\frac{1}{n}\sum_{j=0}^{n} h\left(\frac{j}{n}\right) - \int_{0}^{1} h(x)\,dx$$
$$= \frac{1}{2n}(h(0) + h(1)) + \frac{1}{n}\int_{0}^{1}\left(nx - [nx] - \frac{1}{2}\right)h^{(1)}(x)\,dx.$$

**7. Conclusions.** In this paper we have studied a nonparametric estimator for the pole of a long-memory process under mild conditions on the spectral density $f(\lambda)$. Specifically, we have only assumed that $f(\lambda) \sim C|\lambda - \lambda^0|^{-\alpha}$ with $C > 0$, but smooth elsewhere, and where $\alpha$, the memory parameter, belongs to the interval $(0,1)$. We have shown that the estimator $\hat{\lambda}^0$ of the pole $\lambda^0$ is consistent and we have characterized its limit distribution. More precisely, $\hat{\lambda}^0$, centered around $\lambda^0$ and appropriately renormalized, is asymptotically normal when $\lambda^0 \in (0,\pi)$, whereas if $\lambda^0 = \{0,\pi\}$, the asymptotic distribution is a mixture of a discrete and continuous random variable. In particular, when $\lambda^0 = 0$ the asymptotic distribution takes the value 0 with probability $1/2$ and behaves as a (truncated) normal random variable for positive values. In addition, we have shown that the asymptotic statistical properties of a two-step estimator of $\alpha$ are the same as when $\lambda^0$ is known.

**Acknowledgments.** I thank an Editor, Associate Editor and two referees for their careful reading of the paper and comments that led to a substantially improved version. Also, I have benefited from comments of L. Giraitis. Of course, all remaining errors are my sole responsibility. Finally, I am very grateful to Sue Kirkbride for her patience and excellent typing of this manuscript.

DEPARTMENT OF ECONOMICS
LONDON SCHOOL OF ECONOMICS
HOUGHTON STREET
LONDON WC2A 2AE
UNITED KINGDOM
E-MAIL: f.j.hidalgo@lse.ac.uk